\documentclass[a4paper]{amsart}

\usepackage[UKenglish]{babel} 
\usepackage{amsmath, amsthm, amssymb} 
\usepackage{csquotes, geometry, enumitem} 
\usepackage{xspace, xparse, graphicx} 
\usepackage{tabu} 
\usepackage{tensor}
\usepackage{amsaddr}
\usepackage{geometry}

\theoremstyle{definition}
\newtheorem{defn}{Definition}[section]
\newtheorem{rmk}[defn]{Remark}
\newtheorem{ex}[defn]{Example} 

\theoremstyle{plain}
\newtheorem{thm}[defn]{Theorem}

\newtheorem{prop}[defn]{Proposition}

\newtheorem*{thm*}{Theorem}

\title{Log Floer cohomology for oriented log symplectic surfaces} 
\author{Charlotte Kirchhoff-Lukat} 
\address{Massachusetts Institute of Technology \& KU Leuven\\ 
MIT Department of Mathematics, 77 Massachusetts Ave, Cambridge, MA 02139-4307, USA\\ 
KU Leuven Department of Mathematics, Celestijnenlaan 200B box 2400, 3001 Leuven, Belgium} 
\email{c.kirchhofflukat@kuleuven.be} 

\NewDocumentCommand\party{mmO{}}{\frac{\partial^{#3} #1}{\partial^{#3} #2}}
\NewDocumentCommand\dx{}{\operatorname{d\!}} 
\NewDocumentCommand\res{}{\operatorname{res}}
\NewDocumentCommand\abs{m}{\lvert #1 \rvert} 
\NewDocumentCommand\End{}{\operatorname{End}}

\NewDocumentCommand\Znum{}{\mathbb{Z}}
\NewDocumentCommand\Rnum{}{\mathbb{R}}
\NewDocumentCommand\Cnum{}{\mathbb{C}} 
\NewDocumentCommand\Nnum{}{\mathbb{N}}

\NewDocumentCommand\Imp{}{\operatorname{Im}}
\NewDocumentCommand\Rep{}{\operatorname{Re}}
\NewDocumentCommand\Int{}{\operatorname{Int}}
\NewDocumentCommand\spa{m}{\operatorname{span}\left\{#1\right\}}

\NewDocumentCommand\be{}{\begin{equation}} 
\NewDocumentCommand\ee{}{\end{equation}} 
\NewDocumentCommand\beu{}{\begin{equation*}} 
\NewDocumentCommand\eeu{}{\end{equation*}} 

\NewDocumentCommand\bau{}{\begin{align*}} 
\NewDocumentCommand\eau{}{\end{align*}} 

\NewDocumentCommand\bp{}{\begin{proof}} 
\NewDocumentCommand\ep{}{\end{proof}} 

\NewDocumentCommand\StL{}{\mathbb{D}}

\begin{document} 
\maketitle

\begin{abstract} 
This article provides the first extension of Lagrangian Intersection Floer cohomology to Poisson structures which are almost everywhere symplectic, but degenerate on a lower-dimensional submanifold. 
The main result of the article is the definition of Lagrangian intersection Floer cohomology, referred to as log Floer cohomology, for orientable surfaces equipped with log symplectic structures. We show that this cohomology is invariant under suitable isotopies and that it is isomorphic to the log de Rham cohomology when computed for a single Lagrangian. 
\end{abstract} 

\vspace{0.5cm}
MSC codes: 53D40, 53D17 

Keywords: log symplectic structures, Floer cohomology, real surfaces

Subject classification: symplectic geometry

\tableofcontents

\section{Introduction} 
This work is part of a larger research effort aimed at extending symplectic topology techniques such as Floer theory and Fukaya categories to certain well-behaved Poisson manifolds. Previous work in this area \cite{Pas18,BG22} has taken the approach of integrating the Poisson manifold in question to a symplectic groupoid and studying its symplectic topology. Here, we want to directly define Lagrangian intersection Floer cohomology for a Poisson manifold itself. 

A natural starting point for such an endeavour are manifolds equipped with Poisson structures which are symplectic almost everywhere, but degenerate in a controlled manner on a submanifold. The notion of Lagrangian submanifold extends naturally to this setting, so they are natural candidates for extending the concept of Floer cohomology and Fukaya categories. 

As an initial step, we want to study the setting of real surfaces, where Lagrangian intersection Floer cohomology and the Fukaya category are purely combinatorial. Unfortunately, no stable generalized complex structures (besides those that are just symplectic) exist on surfaces; the smallest dimension with non-trivial examples is 4. But non-trivial log symplectic, also called log Poisson structures, do exist on surfaces; they behave in many ways similarly to certain generalized complex structures on 4-manifolds \cite{Rad02,GMP14}: These Poisson structures are everywhere non-degenerate except on a smooth hypersurface, where they drop in rank by 2, i.e. in the case of a surface, they vanish. 

Consequently, in this work, we take log symplectic structures as proof of concept and as a starting point for extending Floer theory to Poisson manifolds with degeneracies. This article is the first in a pair of research papers: It is a detailed exposition of Lagrangian intersection Floer cohomology for log symplectic surfaces with a thorough introduction to those structures. Since this work constitutes a crossover between Poisson geometry and Floer theory, this article is written to be accessible to readers who are familiar with only one of those subjects. \\
The definition of higher $A_{\infty}$-relations and the construction of the Fukaya category are reserved for the second upcoming article in the series, \emph{The Fukaya category of an orientable log symplectic surface}, \cite{KL23}. 

It is further to be remarked that log symplectic structures exist also on non-orientable surfaces and their theory is as well-understood as in the oriented case \cite{MP18}; the focus of this work on the oriented case is for conciseness and ease of exposition. Every log symplectic structure on a non-orientable surface canonically induces a log symplectic structure on its oriented double cover; using this relationship, the theory can directly be extended to the non-orientable case. 

In Section \ref{sec:olss}, we begin begin by recalling Radko's classification result of log symplectic structures on oriented closed surfaces, as well the refinement using the description of the possible degeneracy loci of these structures first given in \cite{MP18}: Log symplectic structures on a closed surfaces vanish on a collection of embedded disjoint circles $Z=(Z_0,\dots, Z_k)$; the  topological arrangement of these circles up to orientation-preserving diffeomorphism is one invariant of the structure, which can be explicitly described in terms of bipartite graphs. We examine this correspondence in detail and explicitly describe the induced orientation of each vanishing circle in terms of the ambient orientation. 
%

In Section \ref{sec:cohom}, we define the log Lagrangian intersection Floer complex (in the following simply log Floer complex), in particular the \emph{log Floer differential}, for Lagrangians intersecting the vanishing locus $Z$ transversely. We show that this differential does indeed square to zero, and that the resulting cohomology is invariant under Hamiltonian isotopy, with the following crucial particularity: 

In order to ensure transverse intersections of Lagrangians in $Z$ and well-definedness, the \emph{log Floer complex} associated to a pair of Lagrangians $\alpha,\beta$ will always be taken to be $CF(\alpha, \phi^t_H(\beta)$, where $H$ is an \emph{admissible Hamiltonian}. Admissible Hamiltonians $C^{\infty}_+(M)$ constitute a subset of smooth functions with a particular fixed behaviour near the vanishing locus $Z$. This is analogous to Floer theories on non-compact symplectic manifolds such as Wrapped Floer cohomology, where perturbation of one of the Lagrangians with a particular type of Hamiltonian isotopy is also included in the definition of the theory. 

For log symplectic surfaces with multiple symplectic components (i.e. vanishing circles in the interior as opposed to only on a boundary), a crucial feature of this Floer theory is that it includes \emph{crossing lunes}: For two Lagrangians $(\alpha,\beta)$ which have a common intersection point in $Z$ ($\alpha\cap\beta\cap Z \neq \emptyset$), we do not include any pseudoholomorphic discs (here simply called smooth lunes; the holomorphicity condition can be dropped in 2 dimensions) between any intersection points in the symplectic locus and those in $Z$ in the definition of the Floer differential. But we \emph{do include} lunes beginning and ending at intersection points in different symplectic components, which \emph{pass through} a common intersection point in $\alpha\cap\beta\cap Z$. 

\begin{thm*} \textbf{(Theorem \ref{thm:full})}
For $\alpha,\beta\subset (M,\omega,Z)$ closed compact Lagrangians which intersect the degeneracy locus $Z$ transversely and are not null-homotopic under path homotopy leaving $Z$ invariant, the pair $(CF(\alpha,\phi_H(\beta),\partial))$ with $H$ an admissible Lagrangian and $\partial$ the Floer differential defined using crossing lunes, is a cochain complex. 

The associated \emph{log Floer cohomology} is invariant under any deformation of Lagrangians by isotopies that fix the vanishing locus $Z$. 
\end{thm*} 

Crucially, for a single Lagrangian $\alpha$ and a perturbation of $\alpha$ by an admissible Hamiltonian, log Floer cohomology reduces to the log de Rham cohomology of $\alpha$ with respect to $\alpha \cap Z$. This is the natural cohomology theory associated to a manifold with marked hypersurface. 

\begin{thm*}\textbf{\emph{(Proposition \ref{prop:logF})}}
For $\alpha \subset (M,Z)$ a closed embedded Lagrangian in an oriented log symplectic surface that intersects $Z$ transversely in $k$ points, we have 
\begin{align*} 
HF^0(\alpha, \alpha) &\cong \Znum_2 \\ 
HF^1(\alpha, \alpha) &\cong (\Znum_2)^k \oplus \Znum_2,
\end{align*}
while the log de Rham cohomology of $\alpha$ with respect to its endpoints $Z\cap \alpha$ is 
\beu 
H^0(\alpha,\log \alpha \cap Z) \cong \Rnum,\ H^1(\alpha,\log \alpha \cap Z) \cong \Rnum^k \oplus \Rnum
\eeu 
This means that the log Floer cohomology of a single closed Lagrangian $\alpha\cong S^1$ again reduces to the log cohomology of $\alpha$ relative to its intersection with $Z$. 
\end{thm*} 

This constitutes the main argument for why the notion of log Floer cohomology defined here should be considered the correct notion for a Lagrangian intersecting a vanishing circle. 

Finally, in Section \ref{sec:hol}, we describe the Cauchy-Riemann and Floer equations in the setting of log symplectic surfaces and argue how crossing lunes in particular fit into the description as pseudoholomorphic discs. This offers some hints as to how to proceed in higher dimensions, where the theory is no longer combinatorial. 

\paragraph{\textbf{Acknowledgements}} This project has received funding from the European Union’s Horizon 2020 research and innovation programme under the MSCA project \emph{First Steps in Mirror Symmetry for Generalized Complex Geometry} (FuSeGC), grant agreement 887857. 

I would like to thank Paul Seidel for his mathematical mentorship and Abigail Ward and Marco Gualtieri for useful discussions and comments. I further thank the referee for their insightful question, which led to an additional result.

\section{Orientable log symplectic surfaces} \label{sec:olss}
\subsection{Classification} 
Log symplectic structures on closed orientable real surfaces (where all bivectors are automatically Poisson for dimensional reasons) were first systematically studied and classified by O. Radko in \cite{Rad02} under the name \emph{topologically stable Poisson surfaces}: 
\begin{defn} \cite{Rad02} 
A bivector $\pi$ on an orientable surface $M$ defines a \emph{topologically stable Poisson structure} if it vanishes linearly on a finite collection of disjoint embedded circles, and is non-zero everywhere else. 
\end{defn} 
The description as \enquote{topologically stable} refers to the fact that small perturbations of the the bivector will not change the topology of its vanishing locus. 

This class of Poisson structures is dense in the set of all Poisson structures (i.e. all bivectors) on $M$ and provides one of the few examples of a class of Poisson structures that is fully classified. 

Topologically stable Poisson structures on surfaces have since been identified as the 2-dimensional instance of \emph{log Poisson} or \emph{log symplectic structures} (which are shown to be equivalent in \cite{GMP14}), also known as \emph{b-Poisson} or \emph{b-symplectic structures}. 
These structures have been intensively studied in an arbitrary (even) dimension, with the low-dimensional cases (2 and 4) unsurprisingly being best understood. 

\begin{rmk} 
Using the notion of log or b-vector fields and differential forms \cite{Mel93}, the log bivector $\pi$ actually becomes non-degenerate and can be inverted to a log 2-form $\omega:=\pi^{-1}$. As an ordinary 2-form, this is smooth and symplectic away from the vanishing locus of $\pi$. At the vanishing locus of $\pi$, it has a log singularity, hence the name. In \cite{GMP14} the equivalence of log Poisson and log symplectic structures in all dimensions is established. 
\end{rmk} 

\begin{thm} \emph{(\textbf{\cite{Rad02}, Theorem 3})} \label{thm:class}
Topologically stable Poisson structures with $n$ vanishing circles $Z_i$ on a closed oriented surface $M$ are, up to orientation-preserving Poisson isomorphism, completely classified by the following data: 
\begin{enumerate}[label=(\roman*)] 
\item The equivalence class $[Z(\pi)]$ of the set of \emph{oriented} circles $Z(\pi)= \bigsqcup_i Z_i$ up to orientation-preserving diffeomorphism of $M$, 
\item the \emph{modular periods} $\tau_i\in\Rnum_+$ associated to each $Z_i$, 
\item the \emph{regularised Liouville volume} $V(\pi)\in\Rnum$. 
\end{enumerate}
\end{thm} 
\begin{rmk} 
\begin{enumerate} [label=(\roman*)] 
\item In \cite{Rad02}, the author does not specify \emph{which} configurations of circles $\bigsqcup_i Z_i$ can arise as zero loci of topologically stable Poisson structures. This was later done in \cite{MP18} and we are going to detail it in the following subsection. 
\item The topologically stable Poisson bivector $\pi$ induces an orientation on each $Z_i$: The modular vector field of $\pi$ (w.r.t. any volume form on $M$) is tangent to $Z$, and nowhere-vanishing along each $Z_i$. The restriction of the modular vector field along $Z_i$ is independent of the choice of volume form, and in particular orients $TZ_i$. 
\item If $\omega= \pi^{-1}$ is the log symplectic form inverse to the Poisson structure, we can integrate each residue\footnote{The residue of a log symplectic form -- for the purposes of this text simply the inverse of a linearly-vanishing bivector -- is defined as follows: If the vanishing locus $Z_i$ is defined by the linear vanishing of the coordinate $x$, i.e. $Z_i=\{x=0\}$, the vector field $\left.x \party{}{x}\right\vert_{Z_i}$ is independent of the choice of $x$. The residue of $\omega$ along $Z_i$ is the oneform $\res_{Z_i}(\omega)= \iota^*_{Z_i}\left( i_{x\party{}{x}} \omega\right)$.}  $\res_{Z_i}(\omega)$ (which is nowhere vanishing) along $Z_i$ to obtain the \emph{modular period} $\tau_i$. 
\item If $\pi$ is actually nowhere-vanishing, the regularised Liouville volume reduces to the ordinary symplectic volume. For a genuine closed log symplectic surface, the volume of each symplectic component is of course infinite. But since $\omega$ diverges at the same rate, but with opposite sign, on each symplectic component bordering $Z_i$, we can obtain a finite number by regularising the sum of the volumes over all symplectic components of $M$. 
\end{enumerate} 
\end{rmk} 

\subsection{Anatomy} 
From now on, fix a compact orientable surface $M=M^2$, equipped with log Poisson/ log symplectic structure. (Since these notions are equivalent, we will from now on use the terminology interchangeably to describe either the associated bivector or the associated log 2-form.) We denote its bivector by $\pi$ and the corresponding log 2-form by $\omega=\pi^{-1}$. Additionally choose an orientation and volume form $\omega_0$; denote its inverse by $\pi_0=\omega^{-1}$. 
Write $\pi= h \pi_0$; $h$ is a smooth function which vanishes precisely and linearly on $Z= \bigsqcup_i Z_i$, where the $Z_i$ are the disjoint embedded circles that form the degeneracy locus of $\pi$. 
\begin{defn} 
A \emph{symplectic component} of the log symplectic surface $(M,\pi)$ is a connected component of $M\setminus Z$. A symplectic component is called \emph{positive} (or \emph{negative}) if the function $h$ is positive (or negative) on that component. 
\end{defn} 
Obviously the sign of $h$ is fixed on each symplectic component of $M$, and two connected components of $M\setminus Z$ which have a common boundary $Z_i$ have to have opposite signs. This clearly imposes restrictions on those equivalence classes of disjoint circles $\left[Z=\bigsqcup_i Z_i\right]$ in $M$ which are zero loci of some log Poisson structure: 

\begin{prop} \emph{(See \textbf{\cite{MP18}, Theorem 5.5.})}
The possible configurations of vanishing circles $Z_i$ of a log symplectic structure $\pi$ on the oriented surface $(M,\omega_0)$ up to orientation-preserving diffeomorphism are in 1-1 correspondence with bipartite graphs $\Gamma$ whose vertices are decorated with natural numbers $n_j\in \Nnum_0$ s.t. $\sum_j n_j +\operatorname{genus}(\Gamma)=\operatorname{genus}(M)$. 
\end{prop} 
\begin{proof} 
To obtain such a decorated bipartite graph from $(M,\pi)$, assign a vertex to each connected component of $M\setminus Z$. Decorate each vertex with its genus. Finally, draw an edge corresponding to each vanishing circle between the two symplectic components which it bounds. 
The resulting graph is bipartite, because the two symplectic components bounded by each vanishing cycle must have opposite sign. It can be equipped with a vertex 2-colouring by assigning to each vertex the sign of the function $h$ on the corresponding symplectic component. 
Lastly, since $M$ is obtained by gluing all symplectic components together along their boundary, the sum of the genera of the symplectic components and the genus of the graph itself has to be the genus of $M$. 
Clearly, this assignment is invariant under orientation-preserving diffeomorphism. 

Conversely, beginning with a decorated bipartite graph $\Gamma$, construct an oriented surface with a configuration of disjoint circles as follows: For the $j$-th vertex decorated with the natural number $n_j$, take an oriented surface of genus $n_j$ with $k_j$ boundary circles, where $k_j$ is the number of edges connected to the vertex. Then glue these different surfaces with boundary together: Identify those boundary circles that correspond to the same edge in the graph s.t. the orientations match. Keep track of the gluing circle. The result is an oriented surfaces with genus $\sum_j n_j + \operatorname{genus}(\Gamma)$ and a collection of marked circles $\{Z_i\}$. This configuration of circles admits a log Poisson structure whose vanishing circles are precisely the $Z_i$; to find one, choose a vertex 2-colouring of the original graph, then pick any function $h$ that vanishes linearly precisely on the $Z_i$ and whose sign matches the sign of the vertex on each component. Finally, pick a volume form $\omega_0$ matching the orientation and take $\pi = h\omega_0^{-1}$ as the Poisson structure; it will be log Poisson. 
Taking a connected sum is invariant under orientation-preserving diffeomorphism if the orientations of the boundary circles match. 
\end{proof}

As shown in \cite{Rad02} (Section 2.4), any log bivector $\pi$ has the following semilocal form on a cylindrical neighbourhood around each $Z_i$: 
\be 
\pi = c_i x \party{}{x} \wedge \party{}{\theta}, 
\ee 
where $(x\in (-\varepsilon,\varepsilon),\theta\in [0,2\pi))$ are right-handed local coordinates around $Z_i$ (i.e. $Z_i=\{x=0\}$) s.t. $c_i>0$ a positive constant, $\omega_0= \dx x \wedge \dx \theta$, and the modular vector field along $Z_i$ is given by $c_i\party{}{\theta}$, giving the orientation of the vanishing circle $Z_i$. 

\begin{figure}[h]
\centering 
\includegraphics[scale=0.6]{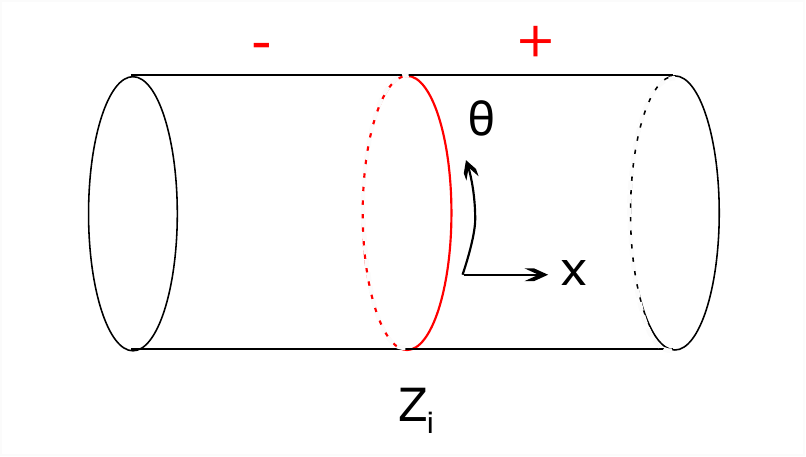} 
\caption{Preferred local coordinates around $Z_i$} 
\end{figure} 

Consider now for a moment $S$, a compact oriented surface with, say, $k$ boundary circles $\{Z_0,\dots,Z_{k-1}\}$ -- think of $S$ as the closure of one of the connected components of $M\setminus Z$ -- equipped with a volume form $\omega_0$. What are possible orientations on the boundary circles that can be induced by a log Poisson structure $\pi$ whose vanishing circles are precisely the $Z_i$? 

\begin{figure}[h] 
\centering 
\includegraphics[scale=0.6]{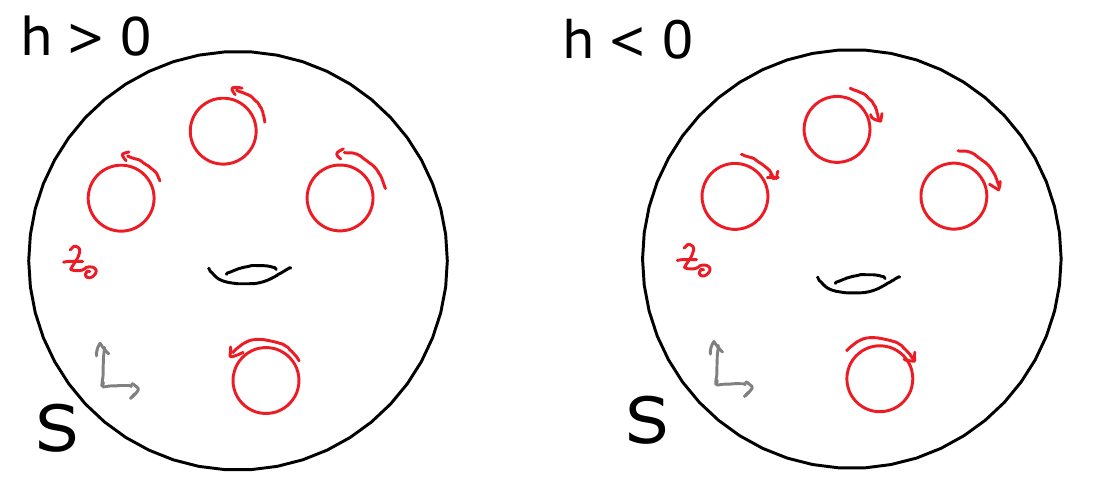}
\caption{With a chosen global orientation of $S$, the sign of the function $h$ determines the orientation of $all$ boundary circles. Here we illustrate the two possible configurations.} \label{fig:bdy_or}
\end{figure} 

\begin{prop} 
If we look at $(S,\partial S=Z=\sqcup_{i=0}^k Z_i)$ with right-handed orientation given by $\omega_0$, any log symplectic structure with vanishing locus $Z$ will have the property that its modular vector field orients the boundary circles in the mathematically positive sense if the orientation of $\pi=h\omega_0^{-1}$ with $h>0$, and in the mathematically negative sense if $h<0$. 
\end{prop} 
\begin{proof} 
Assume that $h>0$; if it is not, flip the global orientation by changing the sign of $\omega_0$. 

Pick two boundary circles, WLOG named $Z_0, Z_1$, and a path $\gamma: [0,1] \rightarrow S$ connecting $Z_0$ to $Z_1$, i.e. $\gamma(0)\in Z_0, \gamma(1)\in Z_1$, such that $\gamma(t) \in \operatorname{Int}(S)$ for all $t\neq 0,1$. Choose a tubular neighbourhood $U$ of the image of $\gamma$. This is diffeomorphic to a strip $[0,1]\times (-\varepsilon, \varepsilon)$. Also pick sets of local coordinates $(x_0,\theta_0)$ and $(x_1,\theta_1)$ as above on neighbourhoods $V_0,V_1$ around both chosen circles. For simplicity, assume that $\gamma(0)$ and $\gamma(1)$ are the respective origins of these coordinates. Then equip $U$ with local coordinates $(a,b)$, which agree with $(x_0,\theta_0)$ on $U\cap V_0$, i.e. 
\beu
x_0=a,\ \theta_0 = b \text{ on } V_0\cap U. 
\eeu
We further require that $\omega_0 = \dx a \wedge \dx b$ on all of $U$. 

\begin{figure}[h] 
\centering 
\includegraphics[scale=0.7]{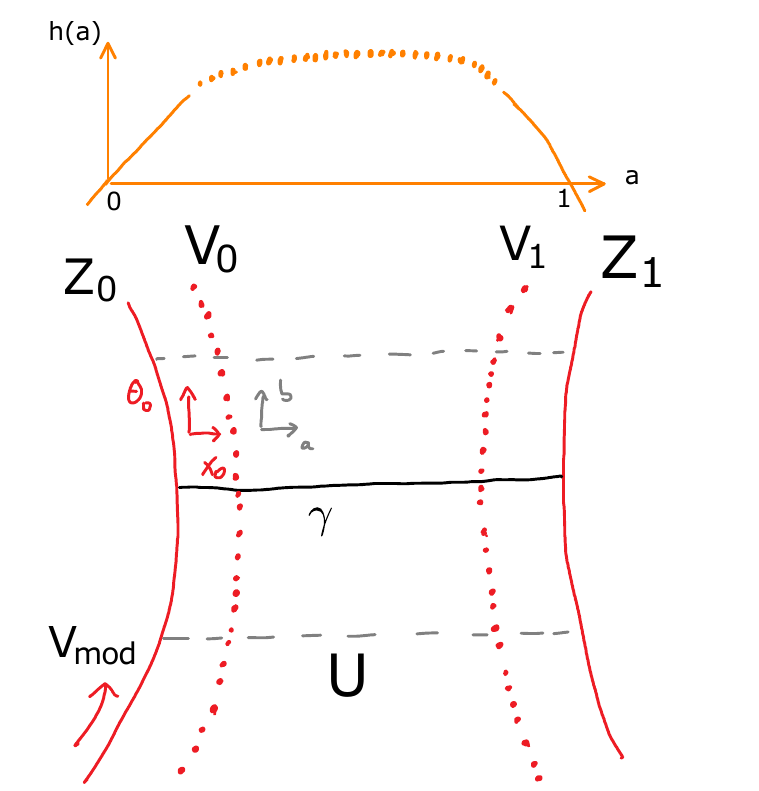}
\caption{Illustration of the path between two boundary circles $Z_0,Z_1$ of $S$ together with chosen coordinate neighbourhoods and a graph of the function $h$ near the endpoints of $\gamma$. The dotted section of the graph indicates that the precise behaviour of $h$ is unknown in this section; we only know that $h>0$.} \label{fig:path1}
\end{figure}

We have started from the assumption that $\pi= h\pi_0$ with $h >0$ on $\operatorname{Int}(S)$. Drawing $\left(\party{}{a},\party{}{b}\right)$ as a right-handed coordinate basis with $\party{}{b}$ pointing upwards, as in Figure \ref{fig:path1}, with the given assumptions, we must have that the modular vector field of $\pi$ restricted to $Z_0$ points upwards as well: The modular vector field on $Z_0$ is given by $c_0\party{}{\theta_0}$ for a positive constant $c_0$, and by assumption $\party{}{\theta}=\party{}{b}$. 

On $V_1$, 
\beu 
\pi = c_1 x_1 \party{}{x_1}\wedge \party{}{\theta_1}\text{ with }c_1>0, \ \pi_0 = \party{}{x_1}\wedge \party{}{\theta_1}. 
\eeu 
At the same time, $\party{h}{a}\vert_{V_1} < 0$, so on the overlap $U\cap V_1$, we must have $a\propto -x_1$. Since $(x_1,\theta_1)$ are by definition also right-handed coordinates, $\party{}{x_1}$ points to the left in our picture Figure \ref{fig:path1} and $\party{}{\theta_1}$ points downwards, i.e. the modular vector field does, too. 

Repeating this process for all boundary circles of $S$ (since $S$ is connected, we can draw paths between $Z_0$ and all other boundary circles, and we can choose these s.t. they are all distinct), this means that when looking at $S$ with the chosen global orientation being right-handed, all boundary circles are oriented in the mathematically positive sense if $h>0$ in the interior of $S$, and correspondingly, that all boundary circles are oriented in the mathematically negative sense if $h<0$. (See Figure \ref{fig:bdy_or} for an illustration.) 
\end{proof} 

For any log symplectic structure, the orientation of all vanishing circles is consequently determined by the overall sign of the bivector $\pi$ -- if we change $\pi$ to $-\pi$, the orientation of all vanishing circles changes, too. Thus, for a given oriented surface $(M,\omega_0)$ with configuration of circles $\{Z_i\}$, the orientation information of a log symplectic structure with vanishing circles $Z_i$ precisely corresponds to the choice of a 2-colouring of the corresponding bipartite graph. 

We can now formulate Radko's classification result more precisely. The following Theorem reformulates Radko's result in \cite{Rad02} explicitly using the graph description pioneered in \cite[Theorem 5.5]{MP18}: 

\begin{thm} \textbf{\emph{(Classification, graph formulation)}} \label{thm:graph}
On a given connected oriented closed surface (with a fixed choice of orientation) $(M,\omega_0)$ of genus $g$, the log symplectic structures $\pi$ are, up to orientation-preserving Poisson isomorphism, in one-to-one correspondence with graphs with the following properties and decorations: 
\begin{enumerate}[label=(\roman*)] 
\item The graph $\Gamma$ is equipped with a 2-colouring. (I.e. it is a bipartite graph.) Additionally, each vertex is decorated with a natural number $g_j\in \Nnum_0$, satisfying $\sum_j g_j + \operatorname{genus}(\Gamma) = g$. 
\item Each edge is labelled with a positive real number $\tau_i$, the modular period of the corresponding vanishing circle. 
\item The entire graph is assigned a real number $V(\pi)$, the regularised Liouville volume. 
\end{enumerate} 
\end{thm} 

\begin{figure}[h] 
\centering 
\includegraphics[scale=0.85]{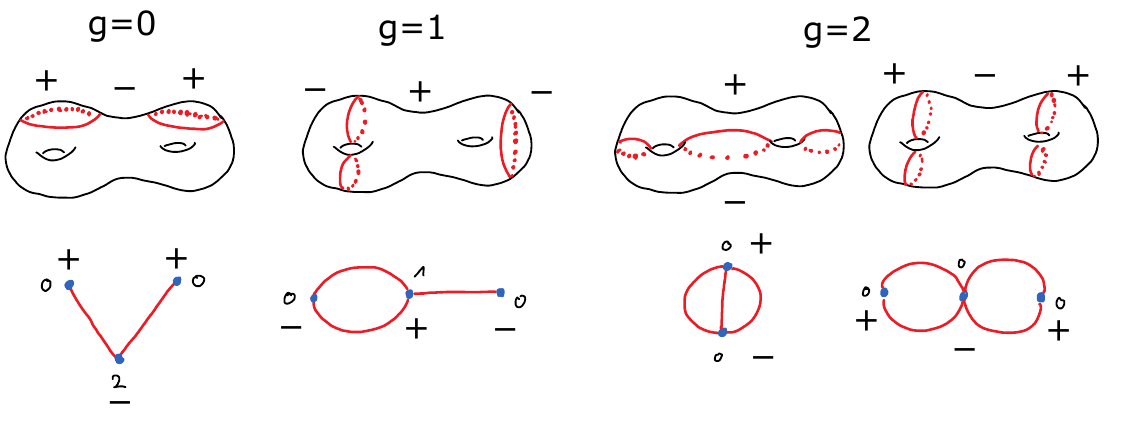} 
\caption{Examples of possible zero loci of log symplectic structures on the orientable surface of genus two, together with the corresponding bipartite graphs. Here, $g$ denotes the genus of the graph. For each example, one of the possible vertex 2-colourings is shown, corresponding to one of the possible relative orientation of the log symplectic structure and global orientation.} \label{fig:genus}
\end{figure} 

\begin{ex}\textbf{The sphere $M=S^2$.} 
As noted in \cite{Rad02}, the log symplectic structures on the sphere can be described in terms of \emph{trees}, i.e. bipartite graphs of genus zero, equipped with a vertex 2-colouring and positive real labels on each edge. There are many log symplectic structures of interest on $S^2$, of increasing complexity, but when we simply say \emph{the log symplectic sphere} or the \emph{standard log symplectic sphere}, we mean the log symplectic structure associated to the following decorated graph: 

\includegraphics[scale=0.4]{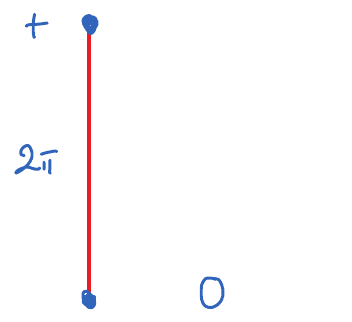}

This is the log symplectic structure with a single vanishing circle along the equator and modular period $\tau=2\pi$ (thus $c=1$). If the vanishing cycle is on the equator, the contributions of the two hemispheres to the regularised Liouville volume precisely cancel, so $V(\pi)=0$. 

Moving the vanishing circle away from the equator, enlarging either the positive or the negative hemisphere, leaves the graph itself invariant, but changes the overall volume to a positive or negative number. 

Adding additional vanishing circles at other level sets of the height function (assuming $S^2$ to be embedded in $\Rnum^3$ with the standard embedding) gives us the \emph{necklace Poisson structures}. Their graphs are trees without any branches. 
\end{ex}

\subsection{Compact Lagrangians and Hamiltonian isotopies} 
The objects of the Fukaya category of a symplectic manifold are (certain classes of) Lagrangians (potentially equipped with additional data). 
Log symplectic structures are genuinely symplectic on a dense subset of the manifold. For our purposes in this article the following naive definition for Lagrangians not entirely contained in the degeneracy locus $Z$ is sufficient:  
\begin{defn} 
A \emph{Lagrangian} submanifold in a log symplectic manifold $(M,Z)$ is a submanifold which is Lagrangian outside $Z$. 
\end{defn} 
We are interested in compact embedded oriented Lagrangians which intersect the degeneracy locus $Z$ transversely. Such Lagrangians have been studied in general dimensions by multiple authors, including in \cite{mythesis,GZ22}.

In 2 dimensions, the Lagrangian condition is of course trivial, so our object of interest inside closed surfaces are embedded oriented circles intersecting $Z$ transversely, in a finite set of points. 
Any Lagrangians that do not intersect $Z$ at all are fully covered by the symplectic theory, so we are going to focus on the case where this intersection is actually non-empty.  
We are also going to consider oriented log symplectic surfaces with boundary, where the vanishing locus of the log Poisson bivector is precisely the boundary -- these are the objects obtained when a closed oriented log symplectic surface is cut apart along all of its vanishing circles. The relevant connected Lagrangians in this setting are embedded closed intervals which intersect the boundary transversely in their endpoints, but which do not intersect it otherwise. (Note that a single circular Lagrangian in a closed log symplectic surface can be separated into several disjoint intervals when the surface is cut apart into components, but in that setting we are still going to identify each connected component as a separate Lagrangian.) 

The notion of Hamiltonian vector field and Hamiltonian isotopy exists for any Poisson structure, and so in particular for any log symplectic structure:The \emph{Hamiltonian vector field} $X_f$ associated to a smooth function $f\in C^{\infty}(M)$ is 
$X_f=\pi(\dx f)$. A \emph{Hamiltonian isotopy} is the isotopy obtained as the flow of a (potentially time-dependent) Hamiltonian vector field. 

Considering again a log symplectic surface $(M^2, Z)$: Recall that in a neighbourhood of each vanishing circle $Z_i$, the log symplectic structure $\pi$ can be written as 
\beu
\pi = c_i x\party{}{x}\wedge \party{}{\theta}. 
\eeu
In these coordinates, the Hamiltonian vector field of $f$ is given by: 
\be 
X_f=\pi(\dx f)= c_i x \party{f}{x}\party{}{\theta} - c_i x \party{f}{\theta}\party{}{x} 
\ee 
Like $\pi$, this of course vanishes linearly on $Z$, and consequently every Hamiltonian isotopy restricts to the identity on $Z$: This means in particular that no two Lagrangians which intersect $Z$ in a different set of points can be Hamiltonian isotopic. Furthermore, any Lagrangian intersecting $Z$ in at least two points (when considering compact embedded circles in a closed surface, this means all Lagrangians that intersect $Z$ at all) is non-contractible under Hamiltonian isotopy. 

Here we can already see some important differences to the purely symplectic case: Lagrangians that intersect $Z$ in different points cannot be related by Hamiltonian isotopy, so log symplectic surfaces have many more equivalence classes of Lagrangians up to Hamiltonian isotopy than symplectic surfaces of the same genus; more Lagrangians need to be considered genuinely different. 

Furthermore, the fact that no compact Lagrangian that intersects $Z$ in at least two points (i.e. whenever an embedded compact Lagrangian intersects $Z$ at all), are contractible under any isotopy fixing $Z$ will turn out to ensure that disc bubbling -- an obstruction to the Floer differential squaring to zero -- cannot arise: 
We address this in Section \ref{sec:cohom}, specifically Remarks \ref{rmk:db2} and \ref{rmk:discbubble}.

\section{Log Floer cohomology for log symplectic surfaces}\label{sec:cohom}
In this section, we mainly follow the exposition of combinatorial Floer cohomology for surfaces in \cite{dSRS14}, which we extend and adapt to the log symplectic setting. We will also make use of the in many ways complementary exposition in \cite{Abo08}. 

Like before, let $(M,Z)$ denote a closed oriented surface (with fixed orientation form $\omega_0$) equipped with log Poisson structure $\pi$ whose vanishing locus is $Z$, and $S$ an oriented surface with $k$ boundary circles $Z=\{Z_0, \dots, Z_{k-1}\}$, to be understood as the closure of one of the symplectic components of $M$. 

As briefly mentioned above, here we focus on Lagrangians whose intersection with $Z$ is non-empty: As soon as one of the compact Lagrangians lies entirely inside a single symplectic component of $M$, the known results on Floer theory for symplectic surfaces apply. 

As in the purely symplectic case, Lagrangian intersection Floer cohomology for log symplectic surfaces is combinatorial, but we include a discussion of the notion of pseudoholomorphic strips and the Cauchy-Riemann equation in the log symplectic setting. 

\subsection{Surface with boundary and a single symplectic component} \label{sec:single}
WLOG assume that the orientations defined by $\omega_0$ and the log symplectic form $\omega=\pi^{-1}$ on $S$ agree, i.e. the sign of the function $h$ s.t. $\pi=h \pi_0$ is positive on the interior of $S$. (The case where the two orientations disagree on a symplectic component only becomes important when there is more than one.) 

Let $(\alpha,\partial \alpha), (\beta,\partial \beta)\subset (S,Z)$ be two compact connected oriented embedded Lagrangians with boundary, i.e. oriented embedded intervals beginning and ending in some boundary component of $S$, which intersect $Z$ transversely. To begin with, assume that 
\be \label{eq:noniso} \alpha\cap\beta\cap Z=\emptyset. \ee
(This implies in particular that $\alpha,\beta$ are non-isotopic under smooth isotopies which are trivial on the boundary.) 

From now on, when referring to smooth isotopy on a log symplectic surface, we always mean those isotopies that restrict to the identity on the vanishing locus. 

\begin{figure}[h] 
\begin{tabu} to \linewidth {X[c]X[c]} 
\includegraphics[scale=0.5]{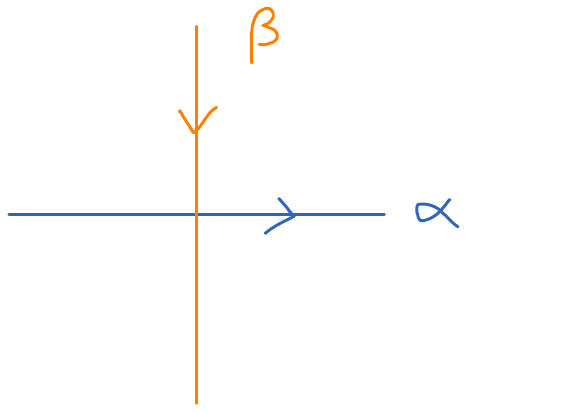} & \includegraphics[scale=0.5]{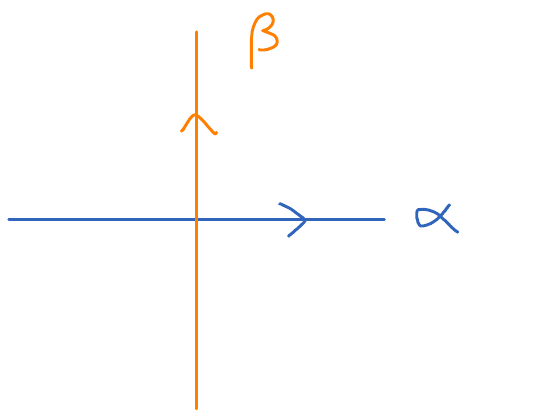} \\ 
\textbf{Case 0:} $\deg(p)=0$ & \textbf{Case 1:} $\deg(p)=1$ 
\end{tabu} 
\caption{Degrees of intersection points of Lagrangians when $h>0$.} \label{fig:deg+}
\end{figure}

\begin{defn} 
The \emph{Floer cochain group} from $\alpha$ to $\beta$ over $\Znum_2$ is 
\be 
CF(\alpha,\beta):= \bigoplus_{p\in \alpha\cap \beta} \Znum_2 p. 
\ee 
Each intersection point $p\in \alpha\cap \beta$ can be assigned a \emph{degree} in $\Znum_2$: 
\be 
\deg(p):= \left\{ \begin{matrix} 
0 &\text{if a neighbourhood of }p\text{ looks like in Figure \ref{fig:deg+}, Case 0}, \\ 
1 &\text{if a neighbourhood of }p\text{ looks like in Figure \ref{fig:deg+}, Case 1}.
\end{matrix}\right.
\ee
This makes $CF(\alpha,\beta)$ into a graded $\Znum_2$-vector space: $CF(\alpha,\beta)= CF^0(\alpha,\beta) \oplus CF^1(\alpha,\beta)$. 
\end{defn} 
\begin{rmk} \begin{enumerate}[label=(\roman*)] 
\item For simplicity, we will restrict the discussion to Floer theory over $\Znum_2$ in this paper; it is possible to work over $\Znum$ instead. 
\item The degree of an intersection point $p$ of course depends on whether $p$ is an element of $CF(\alpha,\beta)$ or $CF(\beta,\alpha)$. 
\item 
Note that since we have assumed $\alpha\cap\beta\cap Z=\emptyset $ (\ref{eq:noniso}), 
all intersection points lie in the interior of $S$. Taking into account the boundary $Z=\partial S$ merely means that the Lagrangians are compact with fixed endpoints under isotopy -- recall that Hamiltonian isotopies of log symplectic surfaces are always trivial on the degeneracy locus $Z$. 
\end{enumerate} 
\end{rmk} 

\begin{figure}[h] 
\centering 
\includegraphics[scale=0.6]{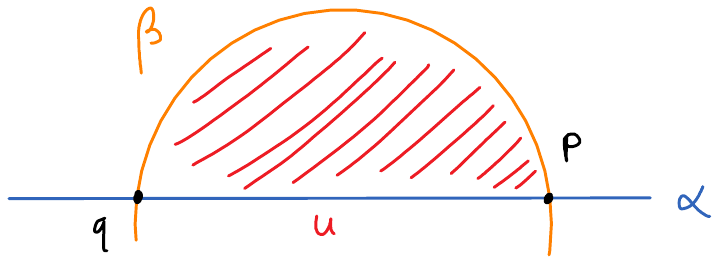}
\caption{A smooth lune from $p$ to $q$ as elements of $CF(\alpha,\beta)$.} 
\end{figure}

\begin{defn} (See Definition 6.1 in \cite{dSRS14}.) 
We call the intersection of the closed unit circle and closed upper half plane $\{z \in \Cnum\vert \abs{z}\leq 1, \Imp(z)\geq 0\}$ in $\Cnum$ the \emph{standard lune} and denote it by $\mathbb{D}$. 

If $p,q\in CF(\alpha,\beta)$, a \emph{smooth lune} from $p$ to $q$ is an equivalence class (up to reparametrisation) of smooth orientation-preserving immersions $u: \StL \rightarrow \Int(S)$ satisfying 
\bau
&u(\StL\cap \Rnum) \subset \alpha,\ &u(\{\abs{z}=1\})\subset \beta \\ 
&u(1)= p,\ &u(-1) =q 
\end{align*}  
and s.t. the corners $p,q$ of the image of $u$ are convex.  
\end{defn} 
Note that smooth lunes always go between intersection points of different degree. 

\begin{defn} 
In this setting, we define the \emph{Floer differential} as follows: 
\be 
\partial: CF(\alpha,\beta) \rightarrow CF(\alpha,\beta),\ \partial(p) = \sum_{q\in \alpha\cap \beta} n(p,q) q, 
\ee 
where $n(p,q)\in \Znum_2$ is the number of smooth lunes from $p$ to $q$ modulo 2. 

The pair $(CF(\alpha,\beta),\partial)$ is then the \emph{Floer cochain complex} or simply \emph{Floer complex} from $\alpha$ to $\beta$. 
\end{defn} 

\begin{thm} \emph{(c.f. Theorems 9.1, 9.2 in \cite{dSRS14}.)} \label{thm:easy}
For non-isotopic $(\alpha,\partial \alpha), (\beta,\partial \beta)$ with $\alpha\cap\beta\cap Z=\emptyset$, the map $\partial$ does indeed make $CF(\alpha,\beta)$ into a differential complex, i.e. $\partial \circ \partial =0$. Furthermore, the associated \emph{Floer cohomology} 
\be 
HF(\alpha,\beta):= \ker(\partial)/\Imp(\partial) 
\ee 
is invariant under isotopies of $S$ which fix $Z$, i.e. 
\beu
HF(\alpha', \beta) = HF(\alpha,\beta),  
\eeu
where $\alpha'$ and $\alpha$ are related by an isotopy that keeps their endpoints fixed. 
\end{thm} 
\bp 
In this particular case, the symplectic argument as outlined in \cite{dSRS14} directly translates in its entirety: $Z$, the boundary and vanishing locus of $\pi$, does not play any role. 
\ep
\begin{rmk}
Theorem \ref{thm:easy} immediately extends to Lagrangians $\alpha=\sqcup_i \alpha_i$ with multiple connected components $\alpha_i$, where each $\alpha_i$ is an embedded interval with endpoints in $Z$ and $\alpha_i\cap \alpha_j=\emptyset$ for $i\neq j$. 
\end{rmk} 

Next, consider the case where the two Lagrangians $\alpha$ and $\beta$ share at least one intersection point with $Z$, i.e. 
\beu
\alpha\cap\beta\cap Z \neq \emptyset. 
\eeu
In this case, issues specific to log symplectic geometry do arise: Should smooth lunes with endpoints in $Z$ be included into the definition of the Floer differential? 
Furthermore, irrespective of whether such lunes are included or not,  the naive definition of the Floer complex as above will no longer yield a cohomology that is invariant under isotopy with fixed endpoints: 

\begin{figure} [h] 
\centering 
\includegraphics[scale=0.5]{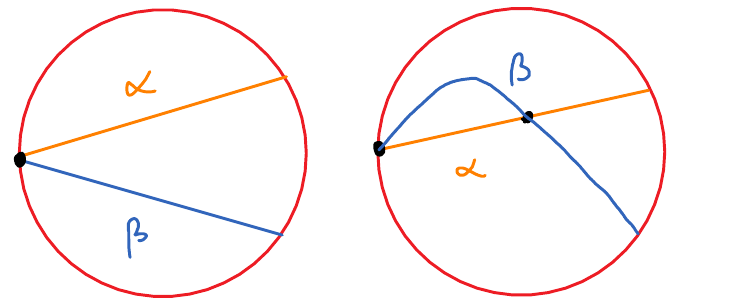} 
\caption{As soon as  $\alpha\cap \beta\cap Z \neq \emptyset$, the Floer complex as naively defined above does no longer result in Floer cohomology that is invariant under isotopy with fixed endpoints.} \label{fig:noninv}
\end{figure} 
\begin{ex} \label{ex:noninv}
Consider $S=D^2$ with log boundary $Z=\partial D^2$. Consider the configurations of Lagrangians $\alpha, \beta$ as in Figure \ref{fig:noninv}, which are related by isotopy with fixed endpoints. On the left, there is only one intersection point, the one inside $Z$, which we call $q$. The Floer differential is trivial, so $HF(\alpha,\beta)=\Znum_2$. On the right, there is an additional intersection point in the interior of $S$, which we denote by $p$. If we define the Floer differential including only lunes which do not touch $Z$, it is again zero, so $HF(\alpha,\beta)=(\Znum_2)^2$. If we do include such lunes, $\partial p = q, \partial q=0$, i.e. $\ker \partial = \Imp \partial$, so $HF(\alpha,\beta)=0$. 
\end{ex} 

\begin{figure}[h] 
\centering 
\includegraphics[scale=0.5]{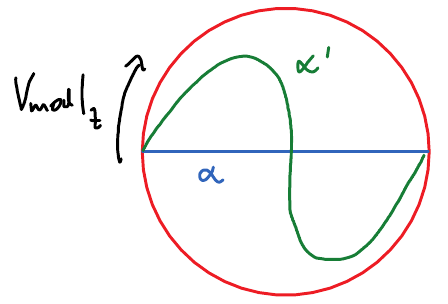} 
\caption{A Lagrangian $\alpha$ is perturbed with an admissible Hamiltonian isotopy to $\alpha'=\phi_H^1(\alpha)$.}\label{fig:admHam}
\end{figure}

In order to fix this, for any $\alpha,\beta$ we set $CF(\alpha,\beta) \equiv CF(\alpha,\phi^t_H(\beta))$, where $\phi^t_H$ is the flow of the Hamiltonian vector field $X_H$ of the function $H\in C^{\infty}_+(S,Z)$. By $C^{\infty}_+(M,Z)$ we mean the subset of smooth functions on $M$ that are of the form 
\be \label{eq:adm}
H(x,\theta)= f(x,\theta) \cdot x, \text{ where } f>0, 
\ee 
on a tubular neighbourhood of $Z$ with coordinates $(x,\theta)$ as above. 
\begin{defn} \label{def:adm} We call functions in $C^{\infty}_+(S,Z)$, i.e. of local form given in equation \ref{eq:adm}, \emph{admissible Hamiltonians}. 
\end{defn} 
We require the time $t$ to be \emph{sufficiently large}, so that further perturbation (i.e. bigger $t$) no longer changes the topological arrangement of $\alpha,\beta$ near $Z$. 
The Hamiltonian isotopy associated to an admissible Hamiltonian will always perturb a Lagrangian in the direction of the modular vector field of $\pi$ on $Z$ near $Z$. This is illustrated in Figure \ref{fig:admHam}. 

If $\alpha\cap\beta\cap Z=\emptyset$, this modification of the definition of Floer complex to include perturbation by an admissible Hamiltonian changes nothing (the original definition already resulted in all isotopies with fixed endpoints inducing quasi-isomorphisms). 
Returning to Example \ref{ex:noninv}, we can now see that $CF(\alpha,\beta)\neq CF(\beta,\alpha)$, but that each corresponding cohomology (whether including the discs reaching $Z$ or not) is invariant under admissible Hamiltonian isotopy. The two configurations in Figure \ref{fig:noninv} now correspond to two different cochain complexes. 
\begin{figure}[h] 
\begin{tabu} to \linewidth {X[c]X[c]} 
\includegraphics[scale=0.5]{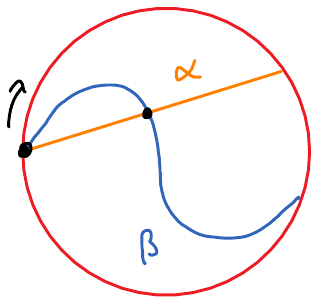} & \includegraphics[scale=0.5]{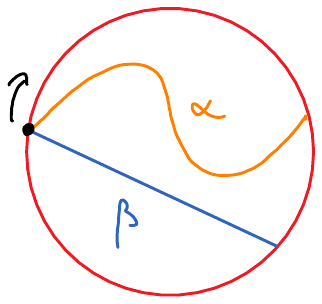} \\ 
$CF(\alpha,\beta)$ & $CF(\beta,\alpha)$ 
\end{tabu} 
\caption{The Floer complexes associated to two Lagrangians in the log symplectic disc with one common intersection point in $Z$.} 
\end{figure} 

This behaviour is nothing new: Various notions of Lagrangian Floer cohomology in non-compact symplectic manifolds perturb Lagrangians going off to infinity using Hamiltonians with specific behaviour near $\infty$. This similarly results non-symmetric Floer complexes and cohomology. Of course $S\setminus Z$ \emph{is} a non-compact symplectic manifold, and Lagrangians that intersect $Z$ go off to infinity in $S \setminus Z$. In this setting of a symplectic surface with log boundary, this is a form of \emph{infinitesimally wrapped} Floer cohomology. 

Note that we can now also unambiguously define the Floer complex of a single Lagrangian in $(S,Z)$: 
\be 
CF(\alpha):= CF(\alpha,\alpha'), 
\ee 
where $\alpha'$ is the image of $\alpha$ under an admissible Hamiltonian isotopy, for example like in Figure \ref{fig:admHam}. 

\begin{defn} \label{defn:logFloer}
The \emph{log Floer differential} is the map
\be 
\partial: CF(\alpha,\phi_H(\beta)) \rightarrow CF(\alpha,\phi^t_H(\beta)), \partial(p) = \sum_{q\in \alpha\cap\phi_H(\beta)} n_{\operatorname{int}}(p,q) q, 
\ee
where $n_{\operatorname{int}}(p,q)\in \Znum_2$ is the number of smooth lunes which lie entirely in $\Int(S)$ from $p$ to $q$ modulo 2 and $\phi_H$ a Hamiltonian isotopy associated to $H\in C^{\infty}_+(S,Z)$, $H=f\cdot x$ with $f$ sufficiently large near $Z$ (in the sense defined above).  
\emph{Log Floer cohomology} is the resulting cohomology $HF(\alpha,\beta)=\ker(\partial)/\Imp(\partial)$. 
\end{defn} 
\begin{rmk} 
In the following, we will simply write $CF(\alpha,\beta)$ for this complex; the modification of the second Lagrangian by an admissible Hamiltonian isotopy is implied. Furthermore, we will simply say \enquote{Floer complex, Floer differential} and \enquote{Floer cohomology} when we mean the log versions defined here. If we want to include lunes with endpoints inside $Z$, we will explicitly point this out. 
\end{rmk} 
Above, we have already assigned a degree to intersection points of $\alpha$ and $\beta$ in the interior of $S$. With the above definition, points in $\alpha\cap\beta\cap Z$ always lie in the kernel of $\partial$. We assign all such points the degree $1$, a decision we will further motivate below. 

\begin{thm} \label{thm:log1}
For two oriented connected Lagrangians $(\alpha,\partial \alpha), (\beta,\partial\beta)$ in the log symplectic surface $(S,Z)$ which intersect $Z$ transversely and precisely in their boundaries, the map $\partial$ as in Definition \ref{defn:logFloer} makes $CF(\alpha,\beta)$ into a differential complex, i.e. $\partial\circ\partial = 0$. Furthermore, the associated log Floer cohomology is invariant under general isotopy leaving $Z$ invariant.
\end{thm} 

\begin{figure}[h]
\centering 
\includegraphics[scale=0.5]{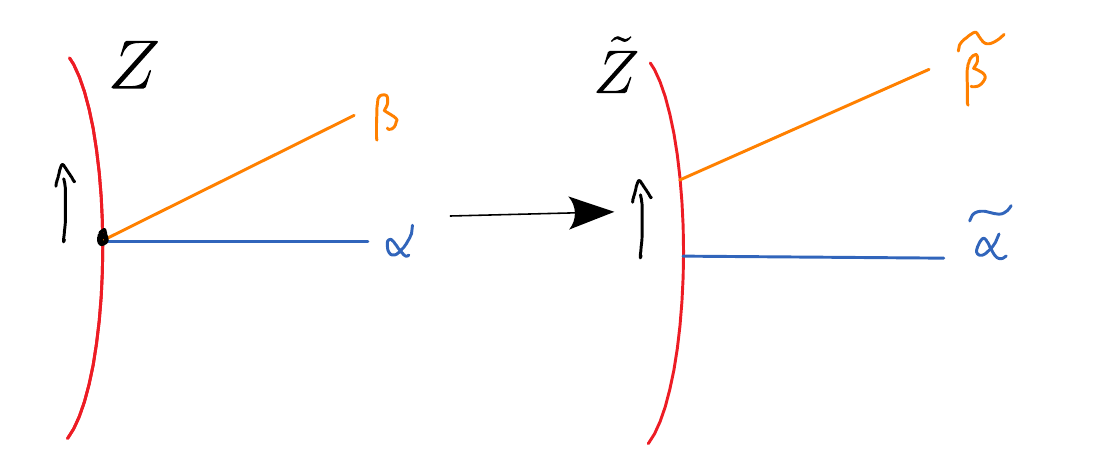}
\caption{On the left: A neighbourhood of $p\in \alpha\cap\beta\cap Z$ when considering $p \in CF(\alpha,\beta)$. On the right, we have performed a real oriented blow-up of $p$, resulting in $\tilde{\alpha}, \tilde{\beta}$, the lifts of the original Lagrangians, intersecting $\tilde{Z}$ in different points.} \label{fig:blow-upbdy} 
\end{figure}

\bp If $\alpha\cap\beta\cap Z= \emptyset$, this reduces to Theorem \ref{thm:easy}, so assume that $\alpha$ and $\beta$ intersect in at least one of their two endpoints, which lie inside $Z$. Let $p$ be such a point. 
According to the definition of $CF(\alpha,\beta)$, perturb $\beta$ sufficiently far using $\phi_H, H\in C^{\infty}_+(M,Z)$ so that $\alpha \pitchfork \beta$. Then a neighbourhood of $p$ looks looks like in Figure \ref{fig:blow-upbdy}. 
We now perform a \emph{real oriented blow-up} of $p$, replacing $p$ by an interval (or circle segment). The resulting surface is still log symplectic ($\omega$ pulls back to a log symplectic form) and $\alpha,\beta$ lift to Lagrangians intersecting the new boundary $\tilde{Z}$ transversely (Furthermore, a real oriented blow-up of a boundary point does not change the topology of $(S,Z)$). However, since $\alpha,\beta$ originally intersected $Z$ in different angles, after blow-up, they intersect $\tilde{Z}$ in different points. \\ 
Of course, $\alpha,\beta\subset S$ are entirely unchanged after the blow-up apart from this, and the log Floer complex $CF(\tilde{\alpha},\tilde{\beta})$ is the same as the original $CF(\alpha,\beta)$. But $\tilde{\alpha},\tilde{\beta}$ satisfy the condition for Theorem \ref{thm:easy}, which we can now apply. 

The cohomology $HF(\alpha,\beta)$ is in fact invariant under general isotopy acting on either Lagrangian, as long as it leaves $Z$ invariant: 
The underlying Floer complex is \emph{always} defined as $CF(\alpha,\phi_H^t(\beta))$, with $t>0$ \emph{sufficiently large}, as defined above (below Definition \ref{def:adm}). If we replace $\beta$ by $\beta'= \psi^{t'}(\beta)$, where $\psi^{t'}$ is not admissible Hamiltonian, then we will in general have to take a bigger $t$ to make the deformation $\phi_H^t(\psi^{t'}(\beta))$ sufficiently large. Applying the blow-up, we find Lagrangians in the setting of Theorem \ref{thm:easy} which are related to the original configuration by isotopy of $\beta$, and consequently have isomorphic Floer cohomology.
If we replace $\alpha$ by $\alpha'=\psi^{t'}(\alpha)$, we note that $CF(\alpha',\beta)\cong CF(\alpha,(\psi^{t'})^{-1}(\beta))$, bringing us back to the situation above.
\ep 

\begin{prop} \label{prop:int}
For $(\alpha, \partial \alpha)\subset (S,Z)$ an embedded interval which transversely intersects $Z$ only in its endpoints, we have
\begin{align*} 
HF^0(\alpha, \alpha) &\cong \Znum_2 \\ 
HF^1(\alpha, \alpha) &\cong (\Znum_2)^2,
\end{align*}
while the log de Rham cohomology of $\alpha$ with respect to its endpoints $Z\cap \alpha$ is 
\beu 
H^0(\alpha,\log \alpha \cap Z) \cong \Rnum,\ H^1(\alpha,\log \alpha \cap Z) \cong \Rnum^2
\eeu 
This means that the log Floer cohomology of a single interval $\alpha$ reduces to the log cohomology of $\alpha$ relative to its intersection with $Z$. 
\end{prop} 
\bp 
According to the Mazzeo-Melrose theorem (Proposition 2.49 in \cite{Mel93}), the log cohomology of a manifold with smooth hypersurface $(M,N)$ is 
\be \label{eq:MaMe}
H^k(M, \log N) \cong H^{k-1}(N) \oplus H^k(M), 
\ee
so for the embedded interval we do indeed obtain: 
\bau 
H^0(\alpha, \log Z\cap \alpha) &= H^0(\alpha) \cong \Rnum_2 \\ 
H^1(\alpha, \log Z\cap \alpha) &= H^0(Z\cap \alpha) \oplus H^1(\alpha) \cong \Rnum^2 
\end{align*} 
By the Weinstein Lagrangian neighbourhood theorem for Lagrangians in log symplectic manifolds (see Theorem 5.18 in \cite{mythesis}), we can take a tubular neighbourhood of $(\alpha,\alpha \cap Z)$ in $(S,Z)$ that is symplectomorphic to a neighbourhood of the zero section in the log cotangent bundle $T^*\alpha(\log \alpha \cap Z)$. Lagrangians near $\alpha$ that are Hamiltonian isotopic to $\alpha$ are precisely those that are graphs of exact log one-forms on $\alpha$. If $x \in [0,1]$ is a coordinate for $\alpha$, the pullback of  $H(x)=-(x-\frac{1}{2})^2+1/4$ is an admissible Hamiltonian. 
\bau
\dx H &= -2\left(x-\frac{1}{2}\right)\dx x \\
&\stackrel{\text{near }x=0}{=} (-2x^2+x)\frac{\dx x}{x} \\ 
&\stackrel{\text{near }x':=1-x=0}{=} (-2(x')^2+x')\frac{\dx x'}{x'} 
\end{align*} 
(We chose the admissible Hamiltonian $H$ s.t. the graph of $\dx H$ in $T^*\alpha(\log \alpha \cap Z)$ is symmetric around the centre of the interval, $x=\frac{1}{2}$.)
Clearly, $\operatorname{graph}(\dx H)\subset T^*\alpha(\log \alpha \cap Z)$ intersects the zero section transversely in $x=0,\frac{1}{2},1$. 
Since there is only one intersection point in the interior of $S$, there are no smooth lunes. Thus 
\bau 
HF^0(\alpha) &\cong \Znum_2 \text{ (Spanned by the interior intersection point.)} \\ 
HF^1(\alpha) &\cong (\Znum_2)^2 \text{ (Spanned by the intersection points at the ends.)}
\end{align*} 
\ep 
This is the motivation for considering a version of Floer cohomology that does not include smooth lunes starting or ending in $Z$ (although we will include lunes passing through $Z$ in the next section): For a single Lagrangian, this log Floer cohomology reduces to log de Rham cohomology. 

\subsection{Closed surface with multiple symplectic components} \label{sec:closed}
In this section we consider a closed log symplectic surface $(M,\pi,Z)$ with non-empty $Z$, i.e. at least 2 symplectic components. Again, fix an orientation $\omega_0$. 

In the previous subsection, we described  the log Floer complex and its cohomology for a symplectic component with positive sign. On those components with negative sign, we have the following changes: 
\begin{figure}[h]
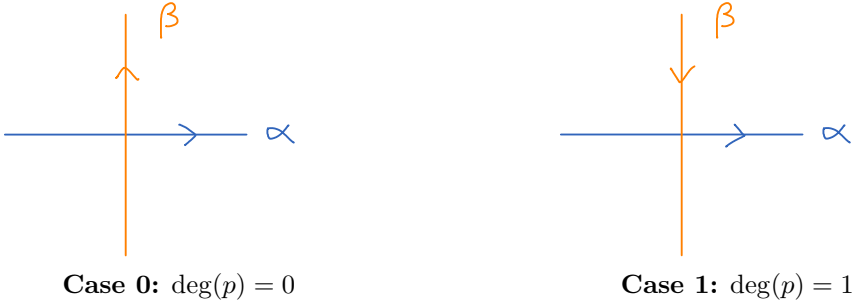
 
\begin{tabu} to \linewidth {X[c]X[c]} 
\includegraphics[scale=0.5]{ori1.pdf} & \includegraphics[scale=0.5]{ori0.pdf} \\ 
\textbf{Case 0:} $\deg(p)=0$ & \textbf{Case 1:} $\deg(p)=1$ 
\end{tabu} 
\caption{Degrees of intersection points of Lagrangians when $h<0$.} \label{fig:deg-}
\end{figure}

\begin{figure}[h] 
\centering 
\includegraphics[scale=0.6]{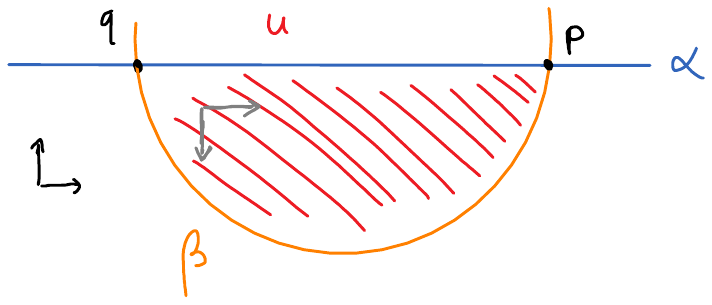}
\caption{A smooth lune from $p$ to $q$ in $CF(\alpha,\beta)$ in a symplectic component where $h<0$. The pairs of basis vectors indicate the orientations of $\StL$ (the standard orientation of $\Cnum$, grey) and the ambient surface (as given by $\omega_0$, black).} \label{fig:neglune}
\end{figure}

\begin{defn} \textbf{(On symplectic components with negative sign)} 
An intersection point in a component of $M\setminus Z$, $p\in CF(\alpha,\beta)$, where $h<0$ is assigned the $\Znum_2$-degree 
\be 
\deg(p):= \left\{ \begin{matrix} 
0 &\text{if a neighbourhood of }p\text{ looks like in Figure \ref{fig:deg-}, Case 0}, \\ 
1 &\text{if a neighbourhood of }p\text{ looks like in Figure \ref{fig:deg-}, Case 1}.
\end{matrix}\right.
\ee

For two intersection points $p,q\in CF(\alpha,\beta)$ in a symplectic component with $h<0$, a \emph{smooth lune} from $p$ to $q$ is an equivalence class (up to reparametrisation) of smooth orientation-reversing\footnote{with respect to the fixed global orientation $\omega_0$} immersions $u: \StL \rightarrow \Int(S)$ satisfying 
\bau
&u(\StL\cap \Rnum) \subset \alpha,\ &u(\{\abs{z}=1\})\subset \beta \\ 
&u(1)= p,\ &u(-1) =q 
\end{align*}  
and s.t. the corners $p,q$ of the image of $u$ are convex. (See Figure \ref{fig:neglune} for an illustration.)  
\end{defn} 
\begin{figure}[h] 
\centering 
\includegraphics[scale=0.6]{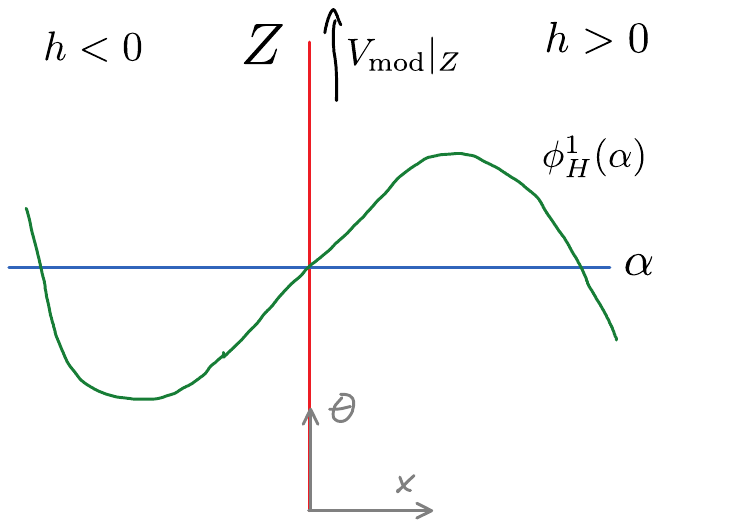}
\caption{Perturbation of a Lagrangian $\alpha$ by an admissible Lagrangian as it crosses $Z$} \label{fig:admHam2}
\end{figure}

Note that we still define $CF(\alpha,\beta)$ to include the sufficiently large perturbation of $\beta$ by an admissible Hamiltonian $H\in C^{\infty}_+(M,Z)$. Figure \ref{fig:admHam2} illustrates how this looks near an an intersection of a Lagrangian with an interior component of $Z$ (which is of course bordered by symplectic components with different signs of $h$).

In addition to these lunes between pairs of intersection points lying in the interior of the same symplectic component, we now also define lunes between intersection points that lie in the interiors of \emph{different} symplectic components. Such lunes only arise in $CF(\alpha,\beta)$ if $\alpha\cap \beta \cap Z \neq \emptyset$. We call such lunes \emph{crossing lunes}. 

Denote by $\StL_k:=(\StL, I_1 \sqcup \dots \sqcup I_k)$ the standard lune with $k\in \Nnum_0$ marked embedded intervals; each $I_j$ is an interval from $\StL\cap \Rnum$ to $\StL\cap \{\abs{z}=s^2+t^2=1\}$ with $s=\text{const}$. 
For any pair $(M,Z)$ of a manifold with embedded hypersurface, denote by $\mathcal{E}_Z$ the restriction of any Euler-like log vector field to the hypersurface $Z$: If $Z$ is given by the vanishing of a coordinate $x$, on any tubular neighbourhood we have have vector field $x\party{}{x}$, which is a nowhere-vanishing section of the log tangent bundle. For the purpose of the definition, pick such an Euler-like vector field both for $(\StL, \sqcup_j I_j)$ and $(M,Z)$. 

\begin{defn} 
For two intersection points $p,q\in CF(\alpha,\beta)$ lying in $M\setminus Z$, a \emph{crossing lune} from $p$ to $q$ is an equivalence class {up to reparametrisation} of smooth maps $u:\StL_k\rightarrow \Int(S)$ (for some $k$) satisfying 
\bau
&u(\StL\cap \Rnum) \subset \alpha,\ &u(\{\abs{z}=1\})\subset \beta \\ 
&u(1)= p,\ &u(-1) =q 
\end{align*}  
and 
\begin{itemize} 
\item $u$ is orientation-preserving wherever $h>0$, and orientation-reversing wherever $h<0$ (this is to say $u$ is orientation-preserving with respect to the orientation of $M\setminus Z$ given by $\omega =\pi^{-1}$), 
\item $(\Imp(u)\cap Z)\subset (\alpha\cap \beta\cap Z)$, 
\item For each $q\in \alpha\cap \beta \cap Z$ that lies in the Image of $u$, $u^{-1}(q)=I_j$ for some $j\in\{1,\dots,k\}$, 
\item $u$ is an immersion on $\StL\setminus u^{-1}(Z)$, and as a map $u_*: T\StL_k(-\log \sqcup_j I_j) \rightarrow TM(-\log Z)$, 
\beu u_*(\mathcal{E}_{\sqcup_j I_j}) = \mathcal{E}_Z ,\eeu
\item and the corners $p,q$ are convex. 
\end{itemize} 
We refer to $\Imp(u)\cap Z$ as the \emph{crossings} of $u$. 
\end{defn} 

Since $u$ is orientation-preserving ($\det(u_*)>0$) on positive symplectic components and orientation-reversing ($\det(u_*)<0$) on negative ones, $\det{u_*\vert_Z}=0$, so $u$ cannot be an immersion on all of $\StL$ if it is to pass through $Z$. 

\begin{rmk}\label{rmk:cl} \begin{enumerate}[label=(\roman*)]
\item Clearly, the above definition reduces to those we have previously given for smooth lunes if $u^{-1}(Z)=\emptyset$. When we talk about the \emph{set of crossing lunes} or \emph{number of crossing lunes} between two points, we usually mean \emph{all} lunes, including those that do not actually pass through $Z$. When we talk about a specific crossing lune, it will be one where $\Imp(u)\cap Z$ is non-trivial, otherwise we will refer to it as a smooth lune. 
\item In fact, with our definition of $CF(\alpha,\beta)$ including the perturbation of $\beta$ by an admissible Hamiltonian, crossing lunes will only ever have one crossing: 
Assume that we have a crossing lune passing through both $q_1, q_2\in\alpha\cap\beta\cap Z$. Consider the submanifold of $M$ enclosed by $\alpha,\beta$ between $q_1,q_2$. This is diffeomorphic to a disc whose boundary is the union of two arc segments $\bar{\alpha}\subset \alpha$ and $\bar{\beta}\subset \beta$. Consequently, $\bar{\alpha}$ and $\bar{\beta}$ are isotopic with fixed endpoints. But this is a contradiction: If the segments $\bar{\alpha},\bar{\beta}$ are isotopic with fixed endpoints in $Z$ and we are considering $CF(\alpha,\beta)\equiv CF(\alpha,\phi_H(\beta))$, we have already perturbed $\beta$ sufficiently far with an admissible Hamiltonian, which will always result in at least one intersection point of $\alpha,\beta$ in the interior of every symplectic component that a crossing lune through either $q_1$ or $q_2$ would begin or end at. \\ 
Consequently, we will take the domain of each crossing lune to be $(\StL, I_1=\{s=0\})$. 
\end{enumerate} 
\end{rmk}

\begin{defn} \label{def:logFloer2}
Given two non-contractible Lagrangians $\alpha, \beta$ in a closed oriented log symplectic surface $(M,\pi,Z)$, the \emph{log Floer differential} is the map 
\be 
\partial: CF(\alpha,\beta)\equiv CF(\alpha,\phi^t_H(\beta)) \rightarrow CF(\alpha,\beta),\ \partial(p) = \sum_{q\in \alpha \cap \phi_H(\beta)} n_c(p,q)q, 
\ee 
where $\phi^t_H, t>0$ is a sufficiently large admissible Hamiltonian isotopy associated to a $H\in C^{\infty}_+(S,Z)$ and $n_c(p,q)$ is the number of crossing lunes from $p$ to $q$ mod 2. (In particular, $n_c(p,q)=0$ if either $p$ or $q$ lies in $Z$.) 
\end{defn} 

\begin{thm} \label{thm:full}
Let $\alpha,\beta \subset (M,Z)$ be two embedded compact Lagrangians in the closed oriented log symplectic surface $(M,Z)$ that intersect $Z$ transversely and non-trivially. Then the map $\partial$ as in Definition \ref{def:logFloer2} makes $CF(\alpha,\beta)\equiv CF(\alpha,\phi_H(\beta))$ into a differential complex, i.e. $\partial \circ \partial =0$. Furthermore, the associated log Floer cohomology is invariant under general isotopy that leaves $Z$ invariant. 
\end{thm} 

\begin{proof} 
If $\alpha \cap \beta \cap Z=\emptyset$, this reduces to Theorem \ref{thm:easy} in each symplectic component of $(M,Z)$, so assume that $\alpha \cap \beta \cap Z$ contains at least one point. If $q\in \alpha \cap \beta \cap Z$, $\partial(q)=0$, so $\partial^2(q)=0$. 

Take $H\in C^{\infty}_+(M)$ an admissible Hamiltonian and $t>0$ sufficiently large, by which we mean that near each $q\in \alpha \cap \beta \cap Z$, $\beta$ lies at a positive angle to $\alpha$ (in the mathematically positive sense), as seen in the first part of Figure \ref{fig:surgery}. 

In order to prove that $\partial^2(p)=0$ for $p \in \alpha\cap \phi_H^t(\beta) \cap (M\setminus Z)$, 
perform the following surgery on $M$: 
\begin{figure}[h] \label{fig:surgery}
\centering
\includegraphics[scale=0.65]{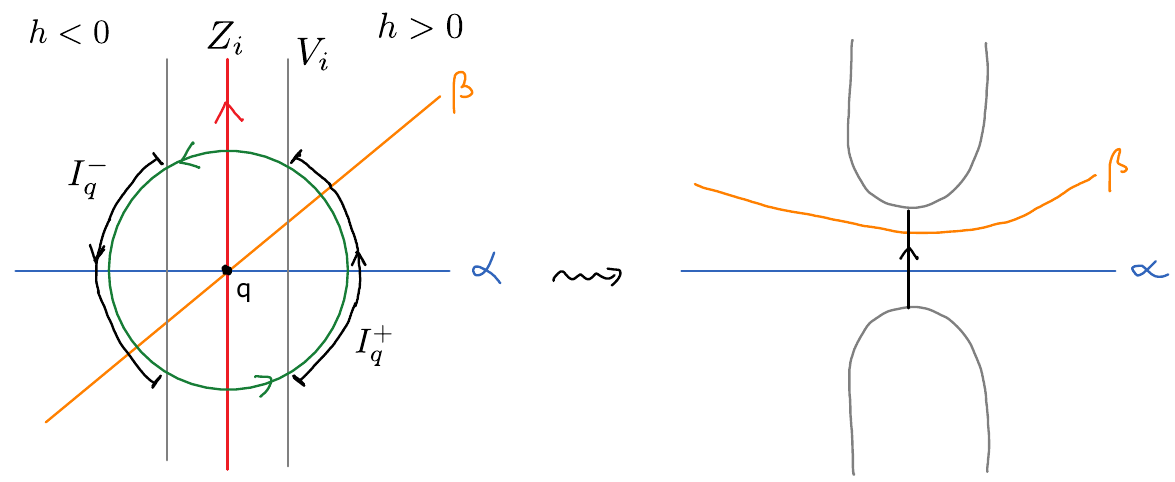}
\caption{This is an illustration of the surgery performed around each $q\in\alpha\cap\beta\cap Z$ identifying $I^+_q$ and $I^-_q$, where on the right we have \enquote{untwisted} the resulting surface. This surgery resolves the intersection point $q$, while still resulting in an orientable surface globally.} \label{fig:surgery}
\end{figure} 

\begin{enumerate}[label=(\roman*)] 
\item Choose a tubular neighbourhood $V_i$ for each of the connected components $Z_i\subset Z$ that contain a $q \in \alpha\cap\beta\cap Z$. Each $V_i$ has to be sufficiently small so as to contain no intersection points of $\alpha$ and $\beta$ except those in $\alpha\cap\beta\cap Z$. 
\item Choose a $D^2$-neighbourhood $U_q$ centred each $q\in \alpha\cap\beta\cap Z$, sufficiently large s.t. 
\beu M':=\left(M\setminus\left(\bigcup_{i} V_i\right)\right)\cup \left(\bigcup_{q\in \alpha\cap\beta\cap Z} U_q\right) \eeu
is connected, but sufficiently small s.t. $U_q$ contains no intersection points of $\alpha,\beta$ other than $q$, nor other points in $(\alpha\cup \beta)\cap Z$. Equip $U_q$ with polar coordinates $(r_q, \vartheta_q)$, chosen s.t. their orientation agrees with that of $\omega$ and $r_q=x_i$. This determines the orientation of the circle $\partial U_q$. 
\item For each $q$, denote by $I^+_q,I^-_q$ the connected components of 
\beu 
(\partial U_q) \cap \left(M\setminus\left(\bigcup_{i} V_i\right)\right), 
\eeu 
where the sign of each is determined by whether $h$ is positve or negative. The $I^{\pm}_q$ inherit an orientation from $\partial U_q$. 
\item Finally, set 
\be \tilde{M} := M\setminus \left(\bigcup_i V_i \cup \bigcup_{q\in\alpha\cap\beta\cap Z} U_q\right)/\sim_{\phi}, 
\ee
where $\phi: \bigcup_q I^+_q \rightarrow \bigcup_q I^-_q$ is an orientation-preserving diffeomorphism that identifies each $I^+_q$ with the corresponding $I^-_q$ such that 
\bau
\phi(I^+_q\cap \alpha)&=I^-_q\cap\alpha \\
\phi(I^+_q\cap \beta)&=I^-_q\cap\beta
\end{align*}
As written, $\tilde{M}$ is a surface with corners, but these can be smoothed out, as can corners that might have arisen in $\alpha$ and $\beta$. 
\end{enumerate} 
This procedure is illustrated in Figure \ref{fig:surgery} around a single point $q\in \alpha\cap \beta\cap Z$. 

Note that $\tilde{M}$ is an oriented surface: We will show this by doing the surgery on each $Z_i$ sequentially, building $\tilde{M}$ step by step, with an oriented surface at each step. Any $Z_i$ separates $M$ into two oriented parts, which we call $M_i^+, M_i^-$ according to the sign of the symplectic component bordering $Z_i$. The boundaries $\partial M^+_i, \partial M^-_i$ inherit orientations from $Z_i$. By the way we have defined the orientations of the $I^{\pm}_q$ (for $q\in Z_i$) we can see that the orientations of the $I^+_q$ agree with the orientation of $\partial M^+_i$, while the orientations of $I^-_q$ are the opposite of $\partial M^-_i$. So now flip the orientations of $M^-_i$ and its boundary, which also changes the sign of every symplectic component of $M^-_i$. When we now glue together the matching $I^\pm_q$ using $\phi$, both the orientations of the glued segments, as well as their normal orientations agree, resulting in a globally oriented surface. This surface is no longer closed, but can be viewed as a surface with boundary. 

It can happen that the resulting $\tilde{\alpha},\tilde{\beta}\subset \tilde{M}$ are no longer circles, but instead collections of intervals ending in $\partial \tilde{M}$, but their intersection points have only changed in that all $q\in \alpha\cap\beta\cap Z$ have been removed. All lunes remain the same, with the exception that lunes with non-trivial crossings are now just smooth lunes: 
\be (CF(\alpha,\beta), \partial) \cong (CF(\tilde{\alpha},\tilde{\beta}), \partial) \oplus \left(\bigoplus_{q\in \alpha\cap\beta\cap Z} \Znum_2 q,\ 0\right) \ee 
We are thus again in the setting of Theorem \ref{thm:easy}. 

Just like for surfaces with a single symplectic component, the cohomology $HF(\alpha,\beta)$ is in fact invariant under general isotopy that leaves $Z$ invariant acting on either Lagrangian. The argument is the same: The underlying Floer complex is \emph{always} defined as $CF(\alpha,\phi_H^t(\beta))$, with $t>0$ sufficiently large; any previous perturbations by arbitrary isotopies fixing $Z$ will not change the surgery. Any changes from the additional isotopy will lie away from the surgery locus; and thus functionally be isotopies within single symplectic components with fixed ends, again allowing us to apply Theorem \ref{thm:easy}. 
\end{proof} 

\begin{rmk} \label{rmk:db2}
This presentation of the proof does not make it explicit, but upon inspection, it becomes clear that no Lagrangian $\alpha$ intersecting $Z$ before the surgery will admit a disc bubble, i.e. a a disc bounded entirely by $\alpha$ -- we have established that pseudoholomorphic discs can only hit $Z$ in discrete points, see also Section \ref{sec:hol}. The surgery is constructed in such a way that the resulting compact Lagrangians in the open honest symplectic surface are non-contractible; there are thus no disc bubbles either. 
\end{rmk} 

The result of Proposition \ref{prop:int} on the log Floer cohomology of a single closed Lagrangian carries over to this setting with multiple symplectic components: 
\begin{prop} \label{prop:logF}
For $\alpha \subset (M,Z)$ a closed embedded Lagrangian in an oriented log symplectic surface that intersects $Z$ transversely in $k$ points, we have 
\begin{align*} 
HF^0(\alpha, \alpha) &\cong \Znum_2 \\ 
HF^1(\alpha, \alpha) &\cong (\Znum_2)^k \oplus \Znum_2,
\end{align*}
while the log de Rham cohomology of $\alpha$ with respect to its endpoints $Z\cap \alpha$ is 
\beu 
H^0(\alpha,\log \alpha \cap Z) \cong \Rnum,\ H^1(\alpha,\log \alpha \cap Z) \cong \Rnum^k \oplus \Rnum
\eeu 
This means that the log Floer cohomology of a single closed Lagrangian $\alpha\cong S^1$ again reduces to the log cohomology of $\alpha$ relative to its intersection with $Z$. 
\end{prop} 

\bp 
Again using Mazzeo-Melrose's result (\ref{eq:MaMe}), we find 
\bau 
H^0(\alpha, \log \alpha \cap Z) &= H^0(\alpha) = H^0(S^1) \cong \Rnum \\ 
H^1(\alpha, \log \alpha \cap Z) &= H^0(Z\cap \alpha) \oplus H^1(\alpha) = H^0(\{\text{pt}\}^{\abs{\alpha \cap Z}})\oplus H^1 (S^1) \cong \Rnum^{\abs{Z\cap \alpha}}\oplus \Rnum
\end{align*} 
Again pick a tubular neighbourhood for $\alpha$ and identify it with an open neighbourhood of the zero section in $T^*\alpha(\log \alpha \cap Z)$. Pick a smooth function $H$ on $\alpha$ with the following properties: 
\begin{enumerate}[label=(\roman*)] 
\item In every tubular neighbourhood $\mathcal{U}_i$ of a $Z_i$ with coordinates $(x_i,\theta_i)$, $H\vert_{\mathcal{U}_i}= x_i$. 
\item On every positive symplectic component of $M$, $H>0$; on every negative symplectic component, $H<0$. 
\end{enumerate} 
Perturbing $\alpha$ with the corresponding Hamiltonian isotopy $\phi^t_{p^*H}$ (where $p: T^*\alpha(\log \alpha \cap Z)\rightarrow \alpha$ is the canonical projection) leads an $\alpha'$ that intersects $\alpha$ in each point in $\alpha\cap Z$, as well as precisely once in every symplectic component that $\alpha, \alpha'$ pass through. If $\abs{\alpha\cap Z}=k$, this means that $CF(\alpha,\alpha')$ is spanned by $q_1,\dots,q_k\in \alpha\cap\alpha'\cap Z$ and $p_1,\dots, p_k \in (\alpha\cap \alpha')\setminus Z$. 
Note that $k$ has to be even. 
WLOG assume that $p_1$ lies in a positive symplectic component, $q_1$ is adjacent, and the numbering proceeds in a mathematically positive sense according to the orientation of $\alpha$. 

If $i\in \{1,\dots,k\}$ is even, two crossing lunes start in $p_i$: One to $p_{i-1}$, one to $p_{i+1}$, where we take $i\in \Znum/k\Znum$. If $i$ is odd, no crossing lunes start in $p_i$. Thus: 
\begin{align*} 
\partial q_i &= 0\ \forall i\\ 
\partial p_i &= \begin{cases} 
p_{i-1} + p_{i+1} &\text{if }i\text{ is odd,} \\ 
0 &\text{if }i\text{ is even.} 
\end{cases} 
\end{align*} 
So we find: 
\bau 
\ker(\partial) &= \spa{q_1,\dots,q_k,p_2,p_4,\dots,p_k, \sum_{j \text{ odd}} p_j} \\ 
\Imp(\partial) &= \spa{p_{i-1}+p_{i+1}\vert i\text{ odd}} \\ 
&= \spa{p_2+p_4, p_4+p_6, \dots, p_{k-2}+p_k} 
\end{align*} 
As a result, the Floer cohomology is 
\beu
HF(\alpha,\alpha') \cong (\Znum_2)^{k+\frac{k}{2}+1}/(\Znum_2)^{\frac{k}{2}-1} \cong \Znum_2^{k+2} \cong H^{\bullet}_{\Znum_2}(\alpha,\log \alpha\cap Z).
\eeu
 
In fact, we can make this isomorphism explicit by associating generators to each other in such a way that it respects degree, i.e. is an isomorphism of graded vector spaces: With the chosen numbering, $p_i$ with odd $i$ are actually intersection points of degree 0, $p_i$ with even $i$ of degree 1. (All $q_i$ are by definition of degree 1). 
We thus obtain an isomorphism that respects degree by assigning 
\be
\begin{aligned} 
 \left[  q_i  \right]&\mapsto [\lambda_i] \\ 
\left[\sum_{i \text{ odd}} p_i\right] &\mapsto [1] \\ 
[p_2] &\mapsto [\zeta], 
\end{aligned} 
\ee
where $\lambda_i$ are the log oneforms locally given by $\lambda_i=\frac{\dx x_i}{x_i}$, $1$ is the constant function and $\zeta$ is the closed nowhere vanishing one-form on the circle.  
\ep

The fact that this version of log Floer cohomology computes log cohomology with respect to the intersection locus with $Z$ motivates our decision not to consider lunes beginning or ending inside $Z$, even though this would also yield a well-defined theory. 
\bigskip


\paragraph{What can we say about the log Floer cohomology of \emph{different} Lagrangians in this setting?} 

Let $\alpha,\beta$ be two Lagrangians with $\alpha\cap Z \neq \emptyset$ and $\beta \cap Z\neq \emptyset$, but $\alpha\cap\beta\cap Z =\emptyset$. Then, just like for ordinary Lagrangian intersection Floer cohomology in a symplectic surface, $HF(\alpha,\beta)$  constitutes an invariant count of minimal intersection points between $\alpha$ and $\beta$, in the sense that its rank over $\Znum_2$ is equal to the minimum number of intersection points under deformations of either Lagrangian by isotopies that leave $Z$ invariant -- after such a deformation, the two Lagrangians will precisely intersect each other so that there are no lunes between intersection points. (See \cite[Chapter 9]{dSRS14}.)

If $\alpha\cap \beta \cap Z \neq \emptyset$, this complicates the situation: 
In order to learn more, it is helpful to consider a pair of Lagrangians which agree near $Z$, by which we mean: For every  $q \in \alpha\cap\beta\cap Z$ there exists an open neighbourhood $U_q\ni q$ such that $\alpha \cap U_q=\beta\cap U_q$. 
This means that for $H\in C^{\infty}_+(M)$ an admissible Hamiltonian, both $CF(\alpha,\phi_H^t(\beta))$ and $CF(\beta, \phi_H^t(\alpha))$ are well-defined -- in the sense of $\phi_H^t$ providing a \emph{sufficiently large} perturbation -- for arbitrarily small $t>0$.

\begin{defn} 
We call a such pair $(\alpha,\beta)$ \emph{adapted}. 
\end{defn} 

If $\alpha,\beta$ are not adapted to begin with, note that we can achieve adaptedness by acting on either $\alpha$ or $\beta$ with suitable isotopies fixing $Z$ which are only non-trivial on a neighbourhood of each $q\in \alpha\cap\beta\cap Z$. 

Now, the sets of generators of $CF(\alpha,\beta)$ and $CF(\beta,\alpha)$ are in general genuinely different in the sense that even with a fixed admissible Hamiltonian $H$ and $t>0$ very small, $\alpha \cap \phi_H^t(\beta)$ and $\phi_H^t(\alpha)\cap \beta \cong \alpha\cap \phi_H^{-t}(\beta)$ can contain a different number of elements for all $t\in (0, \epsilon)$. 

For adapted $(\alpha,\beta,h)$, we distinguish the following types of intersection points: 
\begin{defn} Fix an admissible Hamiltonian $H$. An intersection point $p \in CF(\alpha,\phi_H^t(\beta))$, $t$ small, is called 
\begin{enumerate}[label=(\roman*)]  
\item \emph{inert}: if $p\in \alpha\cap\beta\cap Z$.
\item \emph{ordinary}: if $p$ lies in $M\setminus Z$ and is connected by a path $\gamma: [0,t]\rightarrow M$ in $M$ to a $p_0 \in \alpha\cap \beta\cap (M\setminus Z)$ s.t. $\alpha,\beta$ are transverse at $p_0$, and for all $t'\in [0,t]$, $\gamma(t')$ is an intersection point of $\alpha$ and $\phi_H^{t'}(\beta)$ (i.e. the path $\gamma$ describes the movement of a transverse intersection point of $\alpha$ and $\beta$ as $\beta$ is perturbed by the Hamiltonian isotopy associated to the function $H$).
\item \emph{ephemeral}: otherwise. 
\end{enumerate} 
\end{defn} 

Each \emph{inert} intersection point defines a distinct non-zero class in $HF(\alpha,\beta)$ and $HF(\beta,\alpha)$. 


Each \emph{ordinary} intersection point in $CF(\alpha,\phi_H^t(\beta))$ by definition corresponds to one in $CF(\beta,\phi_H^t(\alpha))$. Of course ordinary intersection points are not necessarily closed, but just like in the situation where $\alpha\cap\beta\cap Z=\emptyset$, we can apply isotopies to $\beta$ to minimise the total number of intersection points of $\alpha$ and $\phi_H^t(\beta)$, while leaving the $U_q$ invariant, so that the remaining ordinary intersection points have no smooth lunes entering or leaving them. Each ordinary intersection point will then be closed and define a distinct non-zero class in $HF(\alpha,\beta)$. 

It is also possible to choose $H$ so that the number of \emph{ephemeral} intersection points is minimised: Near each $q\in \alpha\cap \beta\cap Z$, we then obtain $0,1$ or $2$ ephemeral intersection points directly connected to $q$ through isotopic arcs in both $\alpha$ and $\beta$ without passing through any other intersection points. This number \emph{will} be different for $CF(\alpha,\phi_H^t(\beta))$ and $CF(\beta,\phi_H^t(\alpha))$ and depends on the precise topological arrangement of the two Lagrangians. 

In spite of the differences between $CF(\alpha,\phi_H^t(\beta))$ and $CF(\beta,\phi_H^t(\alpha))$ we now show: In contrast to Lagrangians with boundary inside a symplectic surface with log boundary that we studied in Section \ref{sec:single}, the log Floer cohomology of two such compact Lagrangians in a closed log symplectic surface is always symmetric.

\begin{thm} \label{prop:sym}
Let $(M,\omega,Z)$ a closed oriented log symplectic surface and $(\alpha,\beta)$ embedded closed connected Lagrangians with $\alpha\cap Z \neq \emptyset,\ \beta\cap Z \neq \emptyset$. Then $HF(\alpha,\beta)\cong HF(\beta,\alpha)$. 
\end{thm} 
\begin{proof} 
The only instance where this is still open is when $\beta$ is not isotopic to $\alpha$ under isotopies fixing $Z$, but $\alpha\cap \beta\cap Z\neq \emptyset$.
Without loss of generality, assume that $\alpha,\beta$ are adapted, and fix open neighbourhoods $U_q \ni q$ s.t. $U_q\cap U_{q'}=\emptyset$ for $q\neq q'$. Denote $U:= \bigcup_{q\in \alpha\cap\beta\cap Z} U_q$.

Transform $\beta$ using a Hamiltonian isotopy that leaves the $U_q$ invariant so that in $M\setminus U$, $\alpha,\beta$ intersect transversely in the minimum possible number of intersection points. As established in Theorem \ref{thm:full}, this does not affect log Floer cohomology. 

For each $q\in Z$, further fix an open neighbourhood $V_q \supset U_q$ s.t. $\alpha\cap \beta\cap V_q$ is connected, but $\alpha\cap V_q \neq \beta \cap V_q$. Again make the $V_q$ sufficiently small so that $V_q\cap V_{q'}=\emptyset$ if $q\neq q'$. Fix a tubular neighbourhood $V$ of $Z$ that contains all $V_q$, so that $\alpha\cap \beta\cap V$ has exactly one connected component for each $q\in Z$. 

Then choose an admissible Hamiltonian $H$ that vanishes outside $V$, which on each $V_q$ produces deformations of the type as shown in Figure \ref{fig:vq}, leaving the Lagrangians completely invariant outside the $V_q$. Use this $H$ to define $CF(\alpha,\phi_H^t(\beta)), CF(\beta,\phi_H^t(\alpha))$. 

By construction, all ordinary intersection points in $CF(\alpha,\phi_H^t(\beta))$ and $CF(\beta,\phi_H^t(\alpha))$ agree, only the ephemeral intersection points are different. 

Note: By construction, there are no lunes between ordinary intersection points in this configuration -- if there were, we would not yet have minimised the number of intersection points of $\alpha$ and $\beta$ outside $U$. 
There also cannot be any lunes from an ordinary intersection point to an ephemeral one, or vice versa: If there is such a lune in either $CF(\alpha,\phi_H^t(\beta))$ or  $CF(\beta,\phi_H^t(\alpha))$, this once again means that we have not suitably minimised the number of intersection points between $\alpha$ and $\beta$ outside of $U$; further suitable isotopy leaving $U$ invariant will achieve this. 
Consequently, each ordinary intersection point between $\alpha$ and $\beta$ defines a class in $HF(\alpha,\beta)$ and $HF(\beta,\alpha)$. 

\begin{figure}[t] 
\centering 
\includegraphics[scale=0.5]{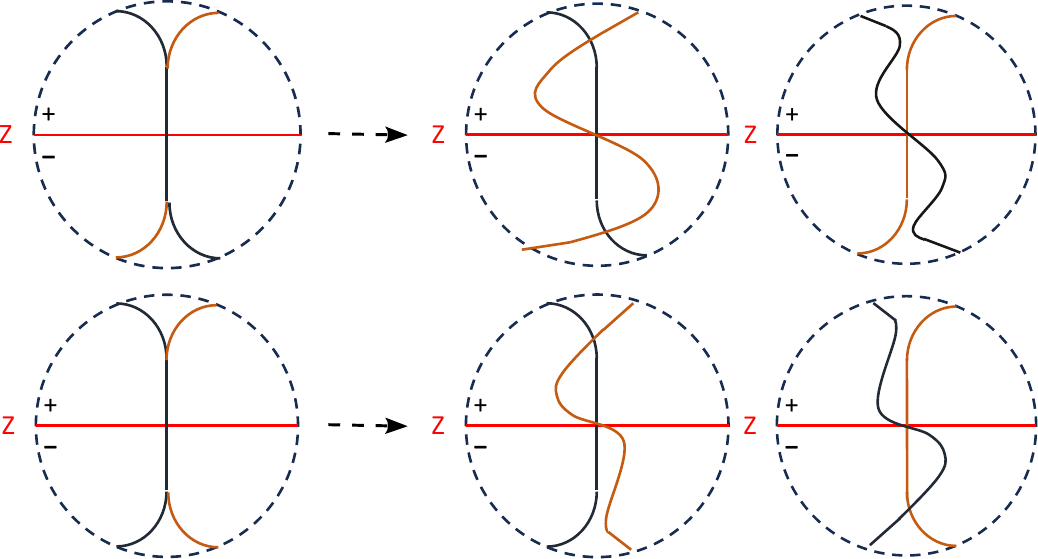}
\caption{An illustration of the ephemeral intersection points arising in the neighbourhoods $V_q$.}  \label{fig:vq}
\end{figure}

If $p_1,p_2 \in CF(\alpha,\phi_H^t(\beta))$ are ephemeral intersection points in the same symplectic component, there cannot be a lune between them by construction.

Now observe, as in Figure \ref{fig:vq}: For each $V_q$, there are four distinct cases: 
\begin{itemize} 
\item \textbf{Case 1:} $\alpha \cap \phi_H^t(\beta)\cap V_q$ contains one ephemeral intersection point $x_{\alpha,\beta}$ in the positive symplectic component. Conversely, $\beta \cap \phi_H^t(\alpha)\cap V_q$ then contains one ephemeral intersection point $y_{\beta,\alpha}$ in the negative symplectic component. 
\item \textbf{Case 2:} $\alpha \cap \phi_H^t(\beta)\cap V_q$ contains two ephemeral intersection points (which are of course connected to $q$ via both Lagrangians): $x$ in the positive symplectic component, $y$ in the negative symplectic component. By construction, there is a crossing lune from $x$ to $y$. \\ 
Correspondingly, in this case $\alpha \cap \phi_H^t(\beta)\cap V_q$ contains \emph{no} ephemeral intersection points.
\item \textbf{Case 3 \& 4:} are obtained by exchanging $\alpha$ and $\beta$ in the above. 
\end{itemize} 
It follows that every ephemeral intersection point $y$ in a negative symplectic component is closed, since no crossing lunes can leave it. 

Search for all maximal sequences $\bar{q}:= (q_1,\dots,q_l)$ of points in $\alpha\cap\beta \cap Z$ s.t. $q_i$ is connected to $q_{i+1}$ by \emph{isotopic arcs} in $\alpha$ and $\beta$ for all $i\in \{1,\dots, l-1\}$. Note that each such maximal sequence has two endpoints, rather than forming a cycle: By assumption $\beta$ is not isotopic to $\alpha$ and both Lagrangians are connected, so there cannot be closed circles where $\alpha,\beta$ are isotopic in each arc. 

The set $\alpha\cap\beta \cap Z$ is a disjoint union of all such $\bar{q}$ and the points $q$ that are not part of any such sequence. 

If $q$ is not part of such a sequence at all, the situation is simplest: If $V_q$ falls into Case 1 or Case 3, the single ephemeral intersection point will define a non-zero cohomology class in both $HF(\alpha,\beta)$ and $HF(\beta,\alpha)$. If $V_q$ falls into Case 2, we have $\partial x=y$, so $[y]=0\in HF(\alpha,\beta)$. No additional non-zero cohomology classes arise in either $HF(\alpha,\beta)$ or $HF(\beta,\alpha)$. (Analogously for Case 4.) 

If $\bar{q}$ is a maximal sequence as described above, the arcs connecting subsequent $q_i\in \bar{q}$ alternately lie in a positive and a negative symplectic component of $M$. There will be exactly one ephemeral intersection point in each such arc in both $CF(\alpha,\phi_H^t(\beta)$ and $CF(\beta,\phi_H^t(\alpha))$.  Denote by $p_i$ the ephemeral intersection point between $q_i$ and $q_{i+1}$. For $2 \leq i\leq l-2$, if $p_i$ is in a positive arc, $\partial p_i= p_{i-1}+p_{i+1}$. As established, all negative ephemeral intersection points are closed by construction, so we find that in both $HF(\alpha,\beta)$ and $HF(\beta,\alpha)$, all $p_i$ in negative arcs define the same cohomology class. 

Next consider $V_{q_1}$: We know that there is an ephemeral intersection point in the arc from $q_1$ to $q_2$ both in $CF(\alpha,\phi_H^t(\beta))$ and in $CF(\beta,\phi_H^t(\alpha))$. There will also be an additional ephemeral intersection point $p_0$ in \emph{either} $CF(\alpha,\phi_H^t(\beta))$ \emph{or} $CF(\beta,\phi_H^t(\alpha))$ on the other side of $q$, where the arcs of $\alpha $ and $\beta$ are no longer isotopic. 
The same applies to $V_{q_l}$, where we denote the resulting ephemeral intersection point by $p_l$. 

In Appendix \ref{app:sym}, we work out the various combinations, where $p_0,p_l$ can each either lie in a negative or a positive symplectic component, and appear in either $CF(\alpha,\phi_H^t(\beta))$ or $CF(\beta,\phi_H^t(\alpha))$. 
In all cases, the maximal sequence $\bar{q}$ contributes either 1 or 0 non-vanishing cohomology classes to both $HF(\alpha,\beta)$ and $HF(\beta,\alpha)$, regardless of how long the sequence $\bar{q}$ is.
\end{proof} 

\begin{rmk} 
In fact, the proof above constitutes a recipe for how to count the rank of general $HF(\alpha,\beta)$ (and, consequently $HF(\beta,\alpha)$). 
\end{rmk}

\begin{ex} 
Consider the arrangement of adapted Lagrangians $\alpha$ and $\beta$ as in Figure \ref{fig:asymex1}. (We can take the picture to show a sufficiently small open subset near any one connected component of $Z$ in any compact log symplectic surface $(M,\omega,Z)$.) 

\begin{figure}[h] 
\centering 
\includegraphics[scale=0.4]{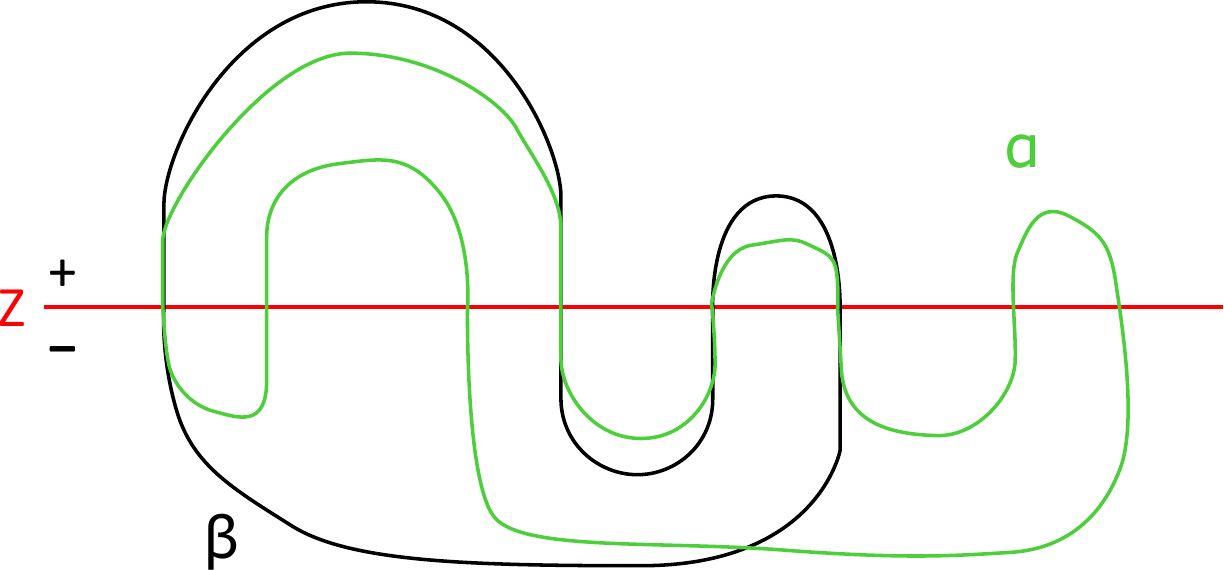}
\caption{Two closed embedded Lagrangians $\alpha, \beta$ in adapted configuration.}\label{fig:asymex1}
\end{figure}

\begin{figure}[h] 
\centering 
\begin{tabu} to \linewidth {X[ll]X[l]} 
(a) & (b) \\
\includegraphics[scale=0.33]{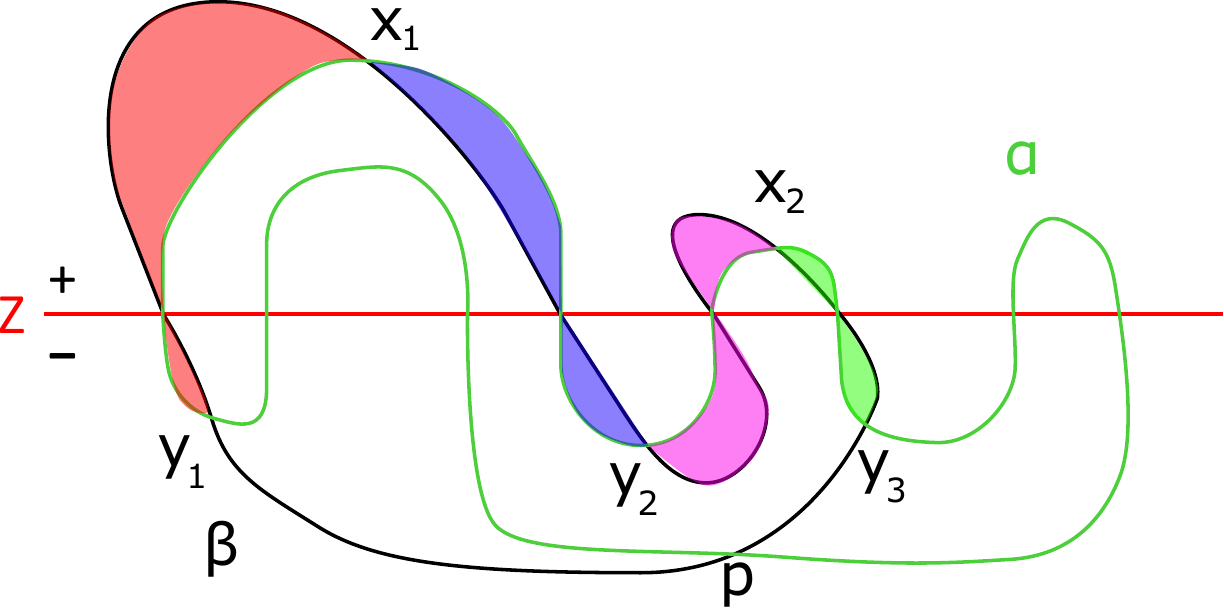} & \includegraphics[scale=0.33]{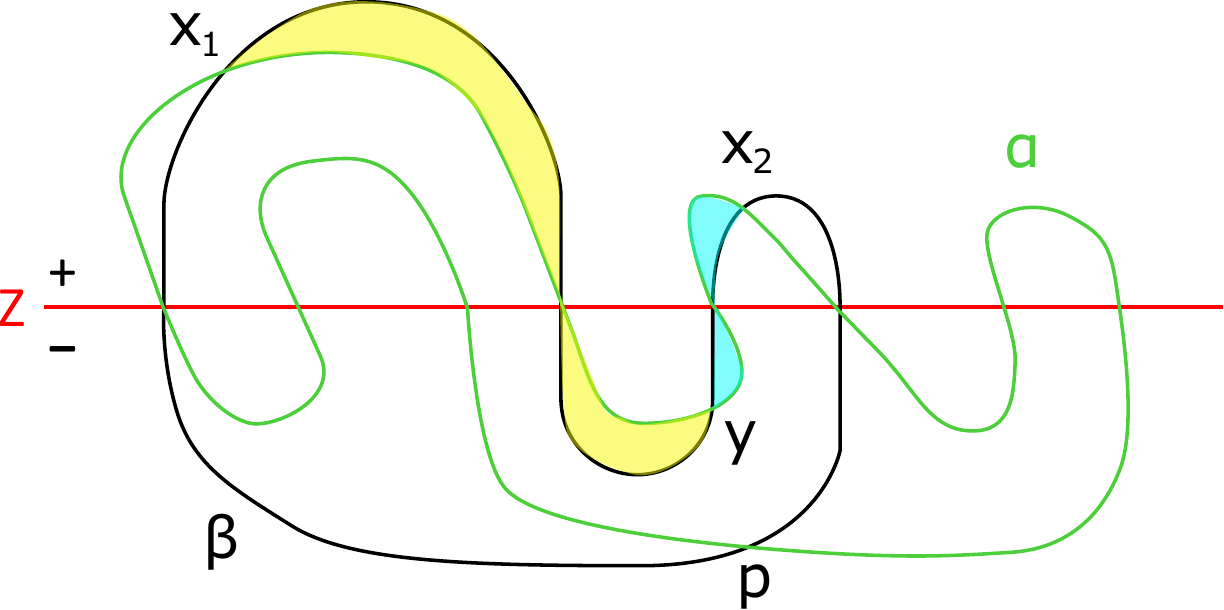}  
\end{tabu} 
\caption{(a)$CF(\alpha,\phi_H^t(\beta))$, (b) $CF(\beta,\phi_H^t(\alpha))$} \label{fig:asymex2} 
\end{figure} 
Figure \ref{fig:asymex2} shows $CF(\alpha,\phi_H^t(\beta))$ and $CF(\beta,\phi_H^t(\alpha))$ for a suitable choice of admissible Hamiltonian $H$. Note that there is a single ordinary intersection point $p$ which defines a non-vanishing class $[p]$ in both $HF(\alpha,\beta)$ and $HF(\beta,\alpha)$. The four inert intersection points $q_1,q_2,q_3,q_4 \in \alpha\cap \beta\cap Z$ are not labelled in the figure. Labelling the $q_i$ as they appear in the figure from left to right, the sequence $\bar{q}=(q_1,\dots, q_4)$ forms a maximal sequence as in the proof above.

We can read off: 
\begin{align*} 
&\text{For }CF(\alpha,\phi_H^t(\beta)): &\text{For }CF(\beta,\phi_H^t(\alpha)): \\ 
&\partial x_1 =y_1 + y_2 &\partial x_1 = y \\ 
& \partial x_2 = y_2 + y_3 &\partial x_2 = y \\ 
& \partial y_1=\partial y_2 =\partial y_3 =0 &\partial y= 0
\end{align*} 
We thus find: 
\begin{align*} 
HF(\alpha,\beta)&= \operatorname{span}([q_1],\dots,[q_4],[p],[y_1]=[y_2]=[y_3])\cong \Znum_2^6 \\ 
HF(\beta,\alpha)&= \operatorname{span}([q_1],\dots,[q_4],[p],[x_1+x_2])\cong \Znum_2^6
\end{align*} 

\end{ex}.

\subsection{A brief aside: Working over the Novikov field} \label{sec:Novikov}
So far in this section, we have always assumed that the Lagrangians under consideration intersect $Z$. Everything we have established so far remains applicable when just one of the Lagrangians does (this is in fact another case where Theorem \ref{thm:easy} applies directly). However, in order for Floer cohomology for pairs of isotopic Lagrangians that do not intersect $Z$ to be defined, we need to work over the Novikov field and modify the definition of the Floer differential to include the symplectic volume of lunes -- this is already the case for symplectic surfaces. Here we briefly outline the necessary modifications: 

\begin{defn} 
The \emph{Novikov field} $\Lambda$ over $\Znum_2$ consists of elements 
\beu 
\sum_{i\in \Znum} a_i T^{\lambda_i}, 
\eeu 
where $a_i\in \Znum_2$ and zero for sufficiently negative $i$, $\lambda_i \in \Rnum$ satisfying $\lim_{i\rightarrow \infty} \lambda_i =\infty$, and $T$ a formal parameter. 
\end{defn} 

The Floer cochain complex over $\lambda$ is $CF(\alpha,\beta)= \bigoplus_{p \in \alpha\cap \phi_H(\beta)} \Lambda\cdot p$, and the log Floer differential is modified to 
\be
\partial(p)=\sum_{\stackrel{q\in \alpha\cap \beta}{u \text{ crossing lune } p\rightarrow q}} T^{\omega(u)} q, 
\ee 
where $\omega(u)= \int_{StL} u^*(\omega)$. Note that we are working in the Novikov field over $\Znum_2$, so two crossing lunes from $p$ to $q$ with equal area will cancel. 

This is only different from the purely symplectic case when we are counting lunes with non-trivial crossing. For these lunes, we need to demonstrate that $\omega(u)$ is well-defined. Furthermore, we need to ensure that the property $\partial \circ \partial =0$ and the invariance under admissible Hamiltonian isotopy are preserved when crossing lunes are involved. 

Recall that the function $h\in C^{\infty}(M)$ s.t. $\pi = h \omega_0^{-1}$ vanishes precisely on $Z$ and does so linearly. Choose a coordinate neighbourhood near $Z_i$ s.t. $\omega =\frac{1}{c_i} \frac{\dx x}{x}\wedge \dx \theta$. Near $q\in \alpha\cap\beta\cap Z_i$, in a potentially smaller coordinate neighbourhood given by $x\in (-\epsilon,\epsilon)$, we can apply admissible Hamiltonian isotopies to $\alpha,\beta$ s.t. 
\bau 
\alpha &= \{(x,0)\} \\ 
\beta &=\{(x,\theta=x)\} 
\end{align*} 
Now assume that a crossing lune passes through $q$. If we remove the tubular neighbourhood $\mathcal{U}$ with $x\in (-\epsilon,\epsilon)$ of $Z_i$, the area of the remaining lune with respect to $\omega$ is clearly finite, as it is the integral of a bounded 2-form over a compact set. We now compute the area of the lune $u$ in $\mathcal{U}$ with respect to $\omega$: 
\bau 
\int_{\mathcal{U}} \omega &= \int_{-\epsilon}^{\epsilon} \int_0^x \frac{1}{x} \dx \theta \dx x \\ 
&= \int_{-\epsilon}^{\epsilon} \dx x \\ 
&= 2\epsilon
\end{align*} 
This shows that crossing lunes have a well-defined area with respect to $\omega$, even though $\omega$ is singular. 

Admissible Hamiltonian isotopies preserve the symplectic area with respect to $\omega$ just like in the symplectic case. What remains is to modify the surgery employed in the proof of Theorem \ref{thm:full} and shown in Figure \ref{fig:surgery} in such a way that the area of crossing lunes does not change in the process: Instead of directly identifying $I_q^{\pm}$, we can insert a strip of the correct width, corresponding to the size of the neighbourhood $U_q$ we removed. Then we can apply the result for the Floer differential over the Novikov field proved in \cite{Abo08}, Lemma 2.11.

\subsection{The Floer equation for log symplectic Floer theory and (pseudo-)holomorphic lunes} \label{sec:hol}
In this section, we motivate the definition of Lagrangian intersection Floer cohomology in terms of the pseudoholomorphic strips satisfying a log Floer equation. 
As a first step, we fix a nice type of compatible almost complex complex structure for the log symplectic form $\omega = \pi^{-1}$: 

\begin{defn} (Compare Definition 4.1.1 in \cite{Alb17}.) An almost complex structure $J\in \End(TM(-\log Z))$ is \emph{cylindrical} if in the adapted coordinate neighbourhood of $Z_i$ with coordinates $(x_i,\theta_i)$ 
\be
J\left(x_i\party{}{x_i}\right)=c_i\party{}{\theta_i},\ J\left(c_i\party{}{\theta_i}\right) = - x_i \party{}{x_i}. 
\ee 
\end{defn} 
Now consider smooth lunes $u: \StL \rightarrow (M,Z)$, where $\StL$ is equipped with the standard complex structure $I$ inherited from $\Cnum$ ($z=s+i t$ is the holomorphic coordinate): 
\beu 
I\left(\party{}{s}\right)= \party{}{t},\ I\left(\party{}{t}\right) = - \party{}{s} 
\eeu
The condition for $u$ to be (pseudo-)holomorphic is 
\beu J \circ u_* = u_* \circ I, \eeu 
which, writing $u$ as $u=(u_x,u_\theta)$ (i.e. suppressing the index $i$), we can rewrite in coordinates as 
\be
\begin{aligned} \label{eq:modCR}
\frac{c}{u_x} \party{u_x}{s}&=\party{u_\theta}{t} \\ 
\frac{c}{u_x} \party{u_x}{t}&=-\party{u_\theta}{s}
\end{aligned} 
\ee
or 
\be
\begin{cases} 
c\party{\log u_x}{s}= \party{u_\theta}{t}; c\party{\log u_x}{t}= - \party{u\theta}{s} &\text{for } u_x >0 \\ 
c\party{\log -u_x}{s}= \party{u_\theta}{t}; c\party{\log -u_x}{t}= - \party{u\theta}{s} &\text{for } u_x <0
\end{cases} 
\ee
which we recognise as the Cauchy-Riemann equations for the complex function $\mathfrak{u}_+:=c\log(u_x) + i u_\theta$, where $u_x>0$ and $\mathfrak{u}_-:= c\log (-u_x) + i u_\theta$ for $u_x<0$.  
By the standard definition of the complex logarithm, 
\bau \mathfrak{u}_+ = \log\left( u_x^c e^{i u_\theta}\right) &=: \log\left(u_+\right) \\ 
\mathfrak{u}_- = \log \left(-u_x^c e^{i u_\theta}\right) &=: \log \left(u_-\right)
\end{align*} 
It thus makes sense to view each connected component of $\mathcal{U}\setminus Z$ (where $\mathcal{U}$ again denotes the standard tubular neighbourhood with coordinates $(x,\theta)$ of $Z$) as equipped with a complex structure equivalent to $\Cnum\setminus \{0\}$ on each connected component, i.e. $z= x^c e^{i \theta}$ where $x>0$ and $z=-x^c e^{i\theta}$ where $x<0$. Alternative holomorphic coordinates for the same complex structures are 
\be \label{eq:holco}
\begin{cases} 
z=x e^{i/c \theta} &\text{for } x >0 \\ 
z= -x e^{i/c \theta}&\text{for } x <0, 
\end{cases} 
\ee 
making $\tau_i=\frac{2\pi}{c_i}$ as the modular period of $Z_i$ more explicit.
The maps $u_{\pm}, \mathfrak{u}_{\pm}$ are holomorphic with respect to these complex structures. 

\begin{figure}[h] 
\centering 
\includegraphics[scale=0.5]{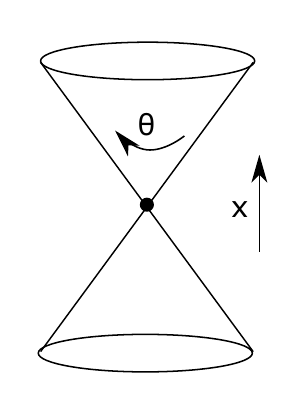} 
\caption{An illustration of the complex structure on $\mathcal{U}\setminus Z$: Each connected component looks like $\Cnum\setminus \{0\}$ in terms of its complex structure; the origin of both components is removed.} \label{fig:complex1} 
\end{figure}

Thus, viewing our original log almost complex structure $J\in \End(TM(-\log Z))$ as a complex structure on $\mathcal{U}\setminus Z$, we obtain something that looks like the neighbourhood of a conical singularity, as in Figure \ref{fig:complex1}. 

In section \ref{sec:Novikov} we demonstrated that crossing lunes have a well-defined, finite area with respect to the log symplectic form $\omega$. They can be viewed as pseudoholomorphic lunes according to (\ref{eq:modCR}) outside $Z$: According to the Riemann mapping theorem, the half discs without endpoints $\StL_-:=\{z\in \mathbb{D}\vert -1<\Rep(z)<0\}$ and $\StL_+:=\{z\in \StL \vert 1>\Rep(z)>0\}$ are (in a neighbourhood of $z=0$, if desired) biholomorphically equivalent to 
\begin{align*} 
\mathbb{S}_- &:= \{s+it\vert s\in (-\infty,\infty),\ 0\leq t \leq e^{-s} \} \\ 
\mathbb{S}_+&:= \{s+it\vert s\in (-\infty,\infty),\ 0\leq t \leq e^{s}\} 
\end{align*} 
respectively.

Recall that we can perturb two Lagrangians $\alpha,\beta$ by admissible Hamiltonian isotopies s.t. inside $\mathcal{U}$, they are given by $\alpha=\{(x,0)\},\ \beta=\{(x,\theta=\frac{1}{c}x)\}$. Thus, when we restrict a crossing lune to each symplectic component of $\mathcal{U}\setminus Z$, we can represent the map $u$ by 
\begin{align*} 
u_-: \mathbb{S}_- &\rightarrow \mathcal{U}_-, (s,t)\mapsto e^{-s+it/c}\\ 
u_+:\mathbb{S}_+ &\rightarrow \mathcal{U}_+, (s,t)\mapsto e^{s+it/c} 
\end{align*} 
respectively, which are holomorphic in the sense of equation (\ref{eq:holco}). 
When factoring these components through $\StL$, we can smooth them out s.t. they extend over $s=0$ to a single smooth map $u$ with $u(s=0)=(x=0,\theta=0)$.

If we further include perturbation of one of the Lagrangians by an admissible Hamiltonian isotopy $H$ in order to ensure $\alpha\pitchfork \beta$ in the setting where this is not already true to begin with, we obtain a \emph{log Floer equation}: 
\begin{align} 
&u_*\left(\party{}{s}\right) + J_t\left(u_*\left(\party{}{t}\right) - X_H \circ u\right) = 0 \\ 
\Leftrightarrow & \label{eq:Floer}
\begin{cases} 
u_x \party{u_{\theta}}{s} + c\party{u_x}{t} + c^2 u_x \left.\party{H}{\theta}\right\vert_u = 0 \\ 
c \party{u_x}{s} - u_x \party{u_{\theta}}{t} + u_x^2 \left. \party{H}{x}\right\vert_u = 0
\end{cases} 
\end{align} 

We thus find that log Floer cohomology is indeed a natural extension of ordinary Floer cohomology for symplectic surfaces, giving us an indication how how to extend it to higher dimensions. 

\begin{rmk} \label{rmk:discbubble}
Note that the equations as written in (\ref{eq:Floer}) explicitly extend across $Z=\{u_x=0\}$. As described above, crossing discs correspond to solutions, and there are \emph{no solutions} $(u_x,u_{\theta})$ to \ref{eq:Floer} that intersect $Z$ in a whole interval, rather than just discrete points. 

This is also the reason why disc bubbling does not arise in the theory as set up here: We consider Lagrangians that are incontractible with respect to path homotopies fixing $Z$. Many of these \emph{are} contractible with respect to general path homotopies on the surface $M$, but the disc bounded by such a Lagrangian will cross $Z$ in at least one interval, and is thus not a pseudoholomorphic disc with respect to the log Floer equation. 
\end{rmk} 

\section{Outlook} \label{sec:outlook}
In  the upcoming second article \cite{KL23}, we define higher $A_{\infty}$-operations for Lagrangians in compact oriented log symplectic surfaces and use, as well as an explicitly constructed collection of Lagrangians, to construct a Fukaya category. We will again argue that this definition constitutes the natural extension of the definition of Fukaya category for a symplectic surface to the log symplectic setting. 

Next, we will attempt to apply the lessons learned from log symplectic surfaces to generalized complex geometry: The aim is to define a \emph{category of generalized complex branes} under certain favourable circumstances. Generalized complex branes constitute the natural submanifolds of generalized complex manifolds \cite{Gua03, Hit03} which extend the notion of \emph{brane} known in symplectic and complex geometry. Since symplectic and complex structures constitute the endpoints of the spectrum of generalized complex structures, it is natural to ask whether homological mirror symmetry as an equivalence of derived categories is fundamentally a generalized complex duality and can be extended to some generalized complex manifolds that are neither symplectic nor complex. However, this requires the notion of a \emph{category of branes} which generalizes the Fukaya category on the symplectic side and the derived category of coherent sheaves on the complex side. 
The first class of manifolds to consider here are compact stable generalized complex 4-manifolds, which in many ways behave similarly to log symplectic surfaces. An additional challenge in this setting comes of course from the fact that the behaviour of pseudoholomorphic curves is much more complicated in 4 dimensions, and the theory will no longer be combinatorial. However, early results indicate that pseudoholomorphic strips between Lagrangians and their moduli spaces still behave well. 

\appendix 
\section{Cases in the proof of Theorem \ref{prop:sym}} \label{app:sym}
Let $\bar{q}=(q_1,\dots,q_l)$ be a sequence of inert intersection points so that each $q_i$ is connected to $q_{i+1}$ by isotopic arcs in $\alpha$ and $\beta$. Denote by $p_i$ the ephemeral intersection point between $q_i$ and $q_{i+1}$, which appears both in $CF(\alpha,\phi_H^t(\beta))$ and $CF(\beta,\phi_H^t(\alpha))$. Denote by $p_0$ and $p_l$ the ephemeral intersection points on the other side of $q_1$ and $q_o$ respectively, which only appear in either $CF(\alpha,\phi_H^t(\beta))$ or $CF(\beta,\phi_H^t(\alpha))$. 

\begin{itemize} 
\item \textbf{(Both $p_1$ and $p_{l-1}$ in negative symplectic component)} In both $CF(\alpha,\phi_H^t(\beta))$ and $CF(\beta,\phi_H^t(\alpha))$, we have: 
\[ \partial p_i = p_{i-1}+p_{i+1} \text{ for }i \text{ even, } \partial p_i =0 \text{ for }i \text{ odd, } i\in \{1,\dots,l-1\}. \] 
 Wlog, assume $p_0 \in CF(\alpha,\phi_H^t(\beta))$, otherwise rename $\alpha,\beta$. Only in $CF(\alpha,\phi_H^t(\beta))$ this gives us the additional identity 
 \[ \partial p_0 = p_1.\] 
	 \begin{enumerate}
	 \item If $p_l \in CF(\alpha,\phi_H^t(\beta))$, this gives us the additional identity $\partial p_l=p_{l-1}$. In this case we find that $[p_i]=0 \in HF(\alpha,\beta)$ for all $i$ odd, but $\partial(p_0+p_2+\dots+p_l)=0$, and $[p_0+p_2+\dots+p_l] \neq 0$. The ephemeral intersection points from the sequence $\bar{q}$ contribute one non-trivial cohomology class in $HF(\alpha,\beta)$.\\ 
	 Now, $CF(\beta,\phi_H^t(\alpha))$ has neither $p_0$ nor $p_l$, so we simply find $[p_1]=[p_3]=\dots=[p_{l-1}]\neq 0 \in HF(\beta, \alpha)$. The ephemeral intersection points from the sequence $\bar{q}$ also contribute one non-trivial cohomology class in $HF(\beta, \alpha)$.
 	\item If $p_l\in CF(\beta,\phi_H^t(\alpha))$, the additional identity $\partial p_l=p_{l-1}$ instead arises there. In this case, we simply find that all closed ephemeral intersection points are also exact in both $CF(\alpha,\phi_H^t(\beta))$ and $CF(\beta,\phi_H^t(\alpha))$, so no additional cohomology class appears. 
 	\end{enumerate}
 
\item \textbf{(Both $p_1$ and $p_{l-1}$ in positive symplectic component)} In both $CF(\alpha,\phi_H^t(\beta))$ and $CF(\beta,\phi_H^t(\alpha))$, we have: 
\[ \partial p_i = p_{i-1}+p_{i+1} \text{ for }i \text{ odd, } \partial p_i =0 \text{ for }i \text{ even, } i\in \{2,\dots,l-2\}. \]
Again, wlog  assume $p_0 \in CF(\alpha,\phi_H^t(\beta))$. Only in $CF(\alpha,\phi_H^t(\beta))$ this gives us the additional identity $\partial p_1=p_0+p_2$, whereas in $CF(\beta,\phi_H^t(\alpha))$: $\partial p_1 = p_2$. 
	\begin{enumerate}
	\item If $p_l \in CF(\alpha,\phi_H^t(\beta))$, this gives us the additional identity $\partial p_{l-1}=p_{l-2}+p_l$, whereas in $CF(\beta,\phi_H^t(\alpha))$ we simply find $\partial p_{l-1}=p_{l-2}$. This means that in $HF(\alpha,\beta)$, the negative ephemeral intersection points define a non-vanishing class, whereas they are trivial in $HF(\beta,\alpha)$. The latter has the non-vanishing class $[p_1+p_3+\dots+p_{l-1}]$ instead. 
	\item If $p_l \in CF(\beta,\phi_H^t(\alpha))$ instead: In both $CF(\alpha,\phi_H^t(\beta))$ and $CF(\beta,\phi_H^t(\alpha))$, the only closed elements coming from $\bar{q}$ are the negative ephemeral intersection points, and they are all exact, so neither $HF(\alpha,\beta)$ nor $HF(\beta,\alpha)$ receives an additional non-vanishing class. 
	\end{enumerate} 
\item \textbf{($p_1$ in a positive symplectic component, $p_{l-1}$ in a negative one)} By renumbering the sequence $\bar{q}$ in the opposite order, this covers the reverse case as well. 
In both $CF(\alpha,\phi_H^t(\beta))$ and $CF(\beta,\phi_H^t(\alpha))$, we have: 
\[ \partial p_i = p_{i-1}+p_{i+1} \text{ for }i \text{ odd, } \partial p_i =0 \text{ for }i \text{ even, } i\in \{2,\dots,l-1\}. \]
Again, wlog  assume $p_0 \in CF(\alpha,\phi_H^t(\beta))$. Only in $CF(\alpha,\phi_H^t(\beta))$ this gives us the additional identity $\partial p_1=p_0+p_2$, whereas in $CF(\beta,\phi_H^t(\alpha))$: $\partial p_1 = p_2$.
	\begin{enumerate}
	 \item If $p_l \in CF(\alpha,\phi_H^t(\beta))$, this gives us the additional identity $\partial p_l=p_{l-1}$. 
	 Consequently the negative ephemeral intersection points are once again the only closed elements in $CF(\alpha,\phi_H^t(\beta))$, and also exact. No additional cohomology class appears. 
	 Now, $CF(\beta,\phi_H^t(\alpha))$ has neither $p_0$ nor $p_l$, but the situation is the same: Only negative ephemeral intersection points are closed, and also exact. 
 	\item If $p_l\in CF(\beta,\phi_H^t(\alpha))$, the additional identity $\partial p_l=p_{l-1}$ instead arises there, in addition to $\partial p_1=p_2$. We again find that $\partial (p_1+p_3+\dots+p_l)=0$, so this element contributes a non-vanishing cohomology class in $HF(\beta,\alpha)$. \\ 
 	In $HF(\alpha,\beta)$ by contrast, the negative ephemeral intersection points $p_i, i$ even, give us the corresponding non-vanishing class. 
 	\end{enumerate}
\end{itemize}  
Thus, in all cases, the maximal sequence $\bar{q}$ contributes either 0 or 1 non-vanishing class to both $HF(\alpha,\beta)$ and $HF(\beta,\alpha)$. 

\bibliographystyle{alpha}
\bibliography{Fukbib.bib}
\end{document}